\documentclass[12pt, leqno]{article}

\usepackage{braket}
\usepackage{hyperref}
\usepackage{theoremref}
\usepackage{amsmath}
\usepackage{amsthm}
\usepackage{amssymb}
\usepackage{amsfonts}
\usepackage{color}
\usepackage[all]{xy} \xyoption{all}
 \usepackage{wrapfig}
\usepackage{enumerate}
\usepackage{enumitem}
\usepackage{slashed}
\usepackage{ytableau}
\usepackage{mathabx}
\usepackage[justification=centering]{caption}

\usepackage{tikz, tikz-3dplot, pgfplots}
\usetikzlibrary{decorations.markings, patterns, decorations.pathmorphing, snakes}
\usetikzlibrary{hobby}

\tikzset{->-/.style={decoration={
markings,
mark=at position #1 with {\arrow{>}}},postaction={decorate}}}

\def\VV{\mathbb{V}}

\def\CC{\mathbb{C}}
\def\DD{\mathbb{D}}

\def\RR{\mathbb{R}}
\def\ZZ{\mathbb{Z}}

\def\TT{{\mathbb{T}}}

\def\Fc{\mathcal{F}}

\def\Nc{\mathcal{N}}

\def\Pc{\mathcal{P}}

\def\Sc{\mathcal{S}}

\def\Uc{{\mathcal {U}}}
\def\Vc{\mathcal{V}}

\def\ba{{\mathbf{a}}}
\def\bb{{\mathbf{b}}}

 \def\bv{{{\bf v}}}

 \def\Aff{{\on{Aff}}}
 
 \def\BM{{\on{BM}}}
 \def\bphi{{\boldsymbol\phi}}
 
 \def\Cor{{\on{Cor}}}
 \def\Crv{{\on{Crv}}}
 
 \def \del{\partial}
 \def\Det{{\on{Det}}}
 \def\Disc{{\on{Disc}}}
 
 \def\Ed{{\on{Ed}}}
 \def\eps{\varepsilon}

 \def\Id{{\on{Id}}}
 \def\Im{{\on{Im}}}
 
\def\Ker{{\on{Ker}}}

\def\Leg{{\on{Leg}}}
\def\lra{\longrightarrow}

\def\on{\operatorname}
\def\oo{{\infty}}
\def\orr{{\on{or}}}

\def\phi{{\varphi}}

\def\PL{{\on{PL}}}
\def\Pr{{\on{Pr}}}

\def\Sk{{\on{Sk}}}

\def\Vert{{\on{Vert}}}
\def\Vol{{\on{Vol}}}

\def\wc{\widecheck}
\def\wt{\widetilde}

\def\0{{\ol{0}}}
\def\1{{\ol{1}}}
\def\(({(\hskip -1mm (}
\def\)){)\hskip -1mm )}
\def\-{{\, \setminus\, }}
\def \= {{\,\, \simeq \,\,}}
 \def\be{\begin{equation}\begin{gathered} }
\def\ee{\end{gathered} \end{equation}}
\def\ed{\end{document}}

\def \beg {\[ \begin{gathered}} 
\def \eng {\end{gathered}\]}

\newtheorem{thm}[equation]{Theorem}
\newtheorem{cor}[equation]{Corollary}

\newtheorem{prop}[equation]{Proposition}

\theoremstyle{definition}

\newtheorem{defi}[equation]{Definition}

\newtheorem{ass}[equation]{Assumption}
\theoremstyle{remark}

\newtheorem{rem}[equation]{Remark}

\newtheorem{ex}[equation]{Example}

\numberwithin{itemcounter}{section}
\numberwithin{equation}{section}
\title{Landau singularities and convex geometry }

\author{Mikhail Kapranov}

\DeclareUnicodeCharacter{2297}{\ifmmode\otimes\else{$\otimes$}\fi}
\DeclareUnicodeCharacter{2194}{\ifmmode\leftrightarrow\else{$\leftrightarrow$}\fi}

\begin{document}
\maketitle
 
 \abstract{We revisit the classic 1959 paper of L. D.  Landau on possible
 singularities of the integral associated  to a Feynman graph.
 Focusing on the real setup, we relate it to convex geometry by representing
 a collection $p$  of incoming momenta as the weighted normals of a convex
 polytope $Q$  via the Minkowski problem. Then,  a regular polyhedral
 subdivision $\Pc$ of $Q$ exhibits $p$ as a Landau singularity labelled
 by the dual graph of $\Pc$ with masses being the areas of the faces. 
 Positivity of the Landau/Feynman/Schwinger multipliers is  interpreted as strict 
 convexity of a PL-function. 
 This gives an interesting class of ``polyhedral'' Landau singularities. 
 
 In the planar case going back to the original 1959 paper, Landau
 graphs can also be identified with plane webs of Gaiotto-Moore-Witten
 that provide a language dual to that of regular polygonal subdivisions.

  }
  
   \tableofcontents  
   
   \vfill\eject

  \section*{Introduction}
   In 1959 L. D. Landau published a seminal paper
   ``On analytic properties of vertex parts in quantum field theory''
   \cite{landau1} where he studied
    possible
 singularities of the integral corresponding to a Feynman graph,
 as a function of the incoming momenta.  
 Since then,
 it became common to refer to all singularities of such Feynman
  integrals
 (and, at times,  even of more general integrals depending on parameters)
  as Landau singularities, see, e.g., \cite{pham}. 
  
  \vskip .2cm
 
 However, the original Landau analysis considered a specific kind of  singularities
 and used certain positivity conditions. 
  In the present paper we 
 relate  these conditions to convex geometry, developing some of
 the ideas present in this analysis. 
 
  \vskip .2cm
 
 First, the classical {\em Minkowski problem}
 \cite{minkowski, klain, gu-yau} allows us to associate
 to any collection $p=(p_i)_{i\in I}$ of incoming momenta (i.e., of  vectors)
 summing up to $0$,  a convex polytope $Q=Q(p)$
 from  which the $p_i$ are recovered as the weighted normals to the
 faces, see \S \ref {sec:minkowski} below.  Any  
 subdivision $\Pc$ of $Q$ into sub-polytopes gives then
 a ``scattering diagram'' $\sigma_\Pc$ labelled  by  $\wt\Gamma$,
 the dual graph of $\Pc$. It is
  obtained by assigning to each edge of $\wt\Gamma$ the 
 weighted normal to the corresponding face $F$ of the corresponding sub-polytope $P$. Such normals sum up to zero for
 each  $P$, ensuring conservation of momentum at each vertex
 of $\wt\Gamma$. 
 
 \vskip .2cm
 
 If the subdivision $\Pc$ is {\em regular} \cite{GKZ}, i.e., admits a strictly
 convex piecewise-linear function $f$, then $\sigma_\Pc$
 represents  the collection of the outer momenta of $\sigma_\Pc$
 as a Landau singularity (Theorem \ref{thm:polyhedral}). The positivity condition for the
 Landau multipliers $\alpha_e$ translates into the convexity
 condition for $f$. The masses of intermediate particles
 are given by  the areas (i.e., volumes in codimension 1)
 of intermediate faces of $\Pc$.  In this way we obtain an interesting
 class of ``polyhedral'' Landau singularities in the real domain. 
 
 \vskip .2cm
 
 This construction represents a higher-dimensional generalization
 of Landau's method of ``schemes'' \cite{landau1, okun}
 associated to  planar Landau graphs (when all momenta lie
 in a $2$-plane).   The interpretation of the set of Landau multipliers as defining a convex piecewise-linear
  function seems to   have  not  been noticed  before.  
   \vskip .2cm
   
   In the planar case, Landau graphs turn out the be the same as
   plane webs of Gaiotto-Moore-Witten \cite{GMW-big, GMW-small}, a concept
   introduced for quite different purposes: study of $\Nc=(2,2)$ supersymmetric
   Landau-Ginzburg theories. The relation between plane webs and
   polygonal subdivisions was explained  earlier in \cite{KKS}. 
   
      \vskip .2cm
   
   In recent years,   the study of scattering  amplitudes
    was connected (partly via
   twistor theory), to various aspects of convex geometry and of its matrix
   generalizations, see  \cite  {arkani-amplitudoh, {AH-EFT}, AH-Henn}. 
   The elementary considerations of the present paper suggest that 
   such connections could possibly  be  traced  to the very foundations
   of quantum field theory. 
   
         \vskip .2cm
         
 The structure of the paper is as follows. In \S \ref {sec:feyn-setup} we recall
         the standard formalism of  integrals associated to Feynman graphs,
         introducing the terminology and notation in terms of the chain complex
         of a graph.  In \S \ref {sec:landau-complex} we review the Landau
         analysis in the complex domain,  introducing the Landau multipliers
         and Landau graphs from direct analysis of the singularity
         conditions. In the next \S \ref {sec:landau-context} we recall some
         other standard features of  Landau's analysis,
          relating them with the so-called Cayley trick in elimination theory  and 
          to the theory of hyperdeterminants \cite{GKZ}. 
          In \S \ref {sec:real-gen} we restrict to the real domain and focus
          on the case when the Landau multipliers are positive, which was the
          original context of \cite{landau1}. In \S \ref {sec:minkowski} we introduce
          the Minkowski problem as a way to encode kinematic data by
          convex polytopes.  After that, \S \ref {sec:RPS-and-RLS}
          introduces a class of {\em polyhedral} Landau singularities
          corresponding to regular subdivisions of the Minkowski polytope. 
          The final \S \ref {sec:planar} is devoted to the planar case when
          the class of polyhedral (polygonal) singularities can be
          characterized directly and when  Landau graphs are
          identified with plane webs of \cite{GMW-big, GMW-small}.

          \vskip .2cm

  This work was supported by  the World Premier International Research Center Initiative (WPI), MEXT, Japan. 


\section{Feynman integrals: the setup}\label{sec:feyn-setup}

In this and  the next section we review  Landau's analysis of singularities of Feynman 
integrals following the original paper \cite{landau1} and
the more systematic treatment in \cite{sato}. 

\paragraph{Background.}  In  quantum  field theory,
the scattering amplitude (the matrix element of the scattering matrix)
is a function  
$S= S(p_1,\cdots, p_n)$ 
of the momenta $p_i$ of incoming and outgoing particles\footnote{For simplicity,
we consider only spinless (scalar) particles.  For particles of nonzero spin,
plane wave states depend not only on the momentum but of additional
discrete parameters.}. 
These $p_i$ are vectors of the dual Minkowski space
$\RR^{3,1}$ satisfy the condition $\sum p_i=0$ (momentum
conservation) and the mass shell conditions $p_i^2=m_i$
being the mass of the $i$th particle. 

\vskip .2cm

In the perturbative approach $S$ is represented as an infinite
sum of contributions
\[
S(p_1,\cdots, p_n) = \sum\nolimits_{\wt\Gamma}\,\,\,  I_{\wt\Gamma}(p_1,\cdots, p_n)
\]
labelled by graphs $\wt \Gamma$ with $n$ ``legs'', edges with one
end free,  labelled by the $p_i$ (the tilde in the notation refers to
the fact that legs are included).

 \begin{figure}[h]
 \centering
 \begin{tikzpicture}[scale=0.6]
 \draw  (0,0) ellipse [x radius=1.5cm, y radius=1cm];
  \fill[gray, opacity=0.2] (0,0) ellipse [x radius=1.5cm, y radius=1cm];
  \node   (1)   at (135: 1.5 and 1) {}; 
    \node   (2)    at (225: 1.5 and 1) {}; 
       \node  (3)   at (-45: 1.5 and 1) {}; 
      \node  (4)   at (0: 1.5 and 1) {}; 
        \node  (5)   at (45: 1.5 and 1) {}; 
        
        \draw (-2.5, 1.2) -- (1.center); 
          \draw (-2.5, -1.2) -- (2.center); 
           \draw (2.5, -1.2) -- (3.center);      
        \draw (3, 0 ) -- (4.center); 
         \draw (2.5, 1.2) -- (5.center);  
         
         \node at (-3, 1.2) {$p_1$}; 
          \node at (-3, -1.3) {$p_2$}; 
           \node at (3, -1.3) {$p_3$};        
             \node at (3.5, 0) {$p_4$}; 
                      \node at (3, 1.2) {$p_5$};    
                      
      \node at (0,0){\large$S$};                   
         \end{tikzpicture}
 \begin{tikzpicture}[scale=0.6]
 \node at (0,1) {\large$=\sum\limits_{\wt\Gamma}
 $
 }; 
 \node at (-2,0){}; 
 \node at (2,0){}; 
 \end{tikzpicture}
  \begin{tikzpicture}[scale=0.6]
  
  \node   (1)   at (135: 1.5 and 1) {\tiny$\bullet$};  
    \node   (2)    at (225: 1.5 and 1) {\tiny$\bullet$}; 
       \node  (3)   at (-45: 1.5 and 1) {\tiny$\bullet$}; 
      \node  (4)   at (0: 1.5 and 1) {\tiny$\bullet$}; 
        \node  (5)   at (45: 1.5 and 1) {\tiny$\bullet$}; 
        
        \draw (-2.5, 1.2) -- (1.center); 
          \draw (-2.5, -1.2) -- (2.center); 
           \draw (2.5, -1.2) -- (3.center);      
        \draw (3, 0 ) -- (4.center); 
         \draw (2.5, 1.2) -- (5.center);  
         
         \node at (-3, 1.2) {$p_1$}; 
          \node at (-3, -1.3) {$p_2$}; 
           \node at (3, -1.3) {$p_3$};        
             \node at (3.5, 0) {$p_4$}; 
                      \node at (3, 1.2) {$p_5$};    
     \draw[line width=1.1](1.center) -- (2.center);      
          \draw [line width=1.1] (2.center) -- (3.center);    
          \draw [line width=1.1] (3.center) -- (4.center);    
      \draw [line width=1.1] (4.center) -- (5.center);     
           \draw [line width=1.1] (1.center) -- (5.center);                 
      \draw [line width=1.1] (2.center) -- (5.center);   
           \draw [line width=1.1] (3.center) -- (5.center);                      
         \end{tikzpicture}
         \caption{An $S$-matrix element and its graph summands}
         \label{fig:S-and-graphs}
 \end{figure}
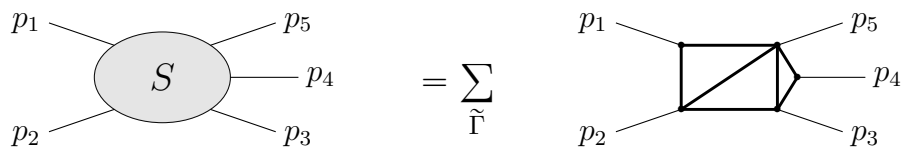
 
  The contribution $I_{\wt\Gamma}$
(the {\em Feynman integral}) is obtained as an integral, over
admissible ways of assigning types of particles (discrete data)
and their
momenta (continuous data) to intermediate edges, of a highly
factorizable function. The latter is built by multiplying:
\begin{itemize}

\item {\em Edge contributions} $G_e$ which are the propagators
(Green's functions) of the corresponding particles.

\item {\em Vertex contributions} $B_v$ obtained from (differential) monomials 
of degree $\geq 3$ (``interaction terms'')
in the Lagrangian, monomials of degree $p$ corresponding to
$p$-valent vertices. 
\end{itemize}

\noindent See below for a more precise formulation. 
Although the integrand is a product, the integral is not since
the integration domain (admissible assignments of internal momenta)
is restricted by requiring conservation of momenta at each vertex. 

\vskip .2cm

Landau \cite{landau1} studied the singularities of $I_{\wt \Gamma}(p_1,
\cdots, p_n)$ as a function of the $p_i$. Note that these are not
the same as singularities of $S(p_1, \cdots, p_n)$ as the
sum over $\wt \Gamma$ is divergent and requires regularization,
even after all the terms are made sense of. 

\paragraph{Formal setup: graphs.}\label{par:graphs}
We consider a finite unoriented graph $\wt\Gamma$,  possibly with legs. For
a formal definition of such objects, see \cite[III.2.3]{manin}. 
We denote by $\Vert(\wt\Gamma)$, $\Ed(\wt\Gamma) \supset L=\Leg(\wt\Gamma)$
the sets of vertices, edges and legs of $\wt\Gamma$, thus thinking of legs as
a particular type of edges (ending at one vertex).  Associating this vertex
to a leg gives a map $\eps: L\to\Vert(\wt\Gamma)$.  
We assume that $\wt\Gamma$ is connected,  has no loops,  that the valence of each
vertex is $\geq 3$ and that $|L|\geq 2$. 

\vskip .2cm

When needed, we will view $\wt\Gamma$ as a $1$-dimensional cell
complex glued out of closed unit intervals $[0,1]$ associated to
all edges which are not legs and of infinite  half-lines $[0, \oo)$
associated to the legs. In particular, each edge of $\wt\Gamma$ has
an affine structure and we can speak about  real valued
functions  on $\wt\Gamma$
which are affine-linear on each edge.  We will call such functions
{\em piecewise-affine}. Similarly for functions with values in any real
vector space. 

\vskip .2cm

We further assume gives a set $\TT$ of ``types of particles''. 
To avoid notational complications, particles are assumed  spinless (scalar).
 Each $t\in\TT$ is assigned a {\em mass} $m_t\in \RR_{\geq 0}$.
We also assume 
that $\wt\Gamma$ to be {\em colored}, 
i.e., equipped with a map $\tau: \Ed(\wt\Gamma) \to \TT$. 
For  $e\in\Ed(\wt\Gamma)$ we write $m_e = m_{\tau(e)}$.  

\vskip .2cm

Let $\Gamma$ be the graph obtained from $\wt\Gamma$ by removing the legs,
so $\Vert(\Gamma)=\Vert(\wt\Gamma)$ and $\Ed(\Gamma) = \Ed(\wt\Gamma) - L$.
Note that $\Gamma$ may have $1$- and $2$-valent vertices. 

 \begin{figure}[h]
 \centering
 \begin{tikzpicture}[scale=0.6]
   \node   (1)   at (135: 1.5 and 1) {\tiny$\bullet$};  
    \node   (2)    at (225: 1.5 and 1) {\tiny$\bullet$}; 
       \node  (3)   at (-45: 1.5 and 1) {\tiny$\bullet$}; 
      \node  (4)   at (0: 1.5 and 1) {\tiny$\bullet$}; 
        \node  (5)   at (45: 1.5 and 1) {\tiny$\bullet$}; 
        
        \draw (-2.5, 1.2) -- (1.center); 
          \draw (-2.5, -1.2) -- (2.center); 
           \draw (2.5, -1.2) -- (3.center);      
        \draw (3, 0 ) -- (4.center); 
         \draw (2.5, 1.2) -- (5.center);  
         
         \node at (-3, 1.2) {}; 
          \node at (-3, -1.3) {}; 
           \node at (3, -1.3) {};        
             \node at (3.5, 0) {}; 
                      \node at (3, 1.2) {};    
     \draw[line width=1.1](1.center) -- (2.center);      
          \draw [line width=1.1] (2.center) -- (3.center);    
          \draw [line width=1.1] (3.center) -- (4.center);    
      \draw [line width=1.1] (4.center) -- (5.center);     
           \draw [line width=1.1] (1.center) -- (5.center);                 
      \draw [line width=1.1] (2.center) -- (5.center);   
           \draw [line width=1.1] (3.center) -- (5.center);   
           
           \node at (0, -2) {\large$\wt\Gamma$}; 
 
 \end{tikzpicture}
  \begin{tikzpicture}[scale=0.6]

 \node at (0,0.2){\huge$\rightsquigarrow$}; 
 
\node at (0, -2.0){};

 \end{tikzpicture}
  \begin{tikzpicture}[scale=0.6]
   \node   (1)   at (135: 1.5 and 1) {\tiny$\bullet$};  
    \node   (2)    at (225: 1.5 and 1) {\tiny$\bullet$}; 
       \node  (3)   at (-45: 1.5 and 1) {\tiny$\bullet$}; 
      \node  (4)   at (0: 1.5 and 1) {\tiny$\bullet$}; 
        \node  (5)   at (45: 1.5 and 1) {\tiny$\bullet$}; 

         \node at (-3, 1.2) {}; 
          \node at (-3, -1.3) {}; 
           \node at (3, -1.3) {};        
             \node at (3.5, 0) {}; 
                      \node at (3, 1.2) {};    
     \draw[line width=1.1](1.center) -- (2.center);      
          \draw [line width=1.1] (2.center) -- (3.center);    
          \draw [line width=1.1] (3.center) -- (4.center);    
      \draw [line width=1.1] (4.center) -- (5.center);     
           \draw [line width=1.1] (1.center) -- (5.center);                 
      \draw [line width=1.1] (2.center) -- (5.center);   
           \draw [line width=1.1] (3.center) -- (5.center);   
           
           \node at (0, -2) {\large$\Gamma$}; 
 
 \end{tikzpicture}
 \caption {Removing legs from a scattering graph.} 
 \label{}
 \end{figure}
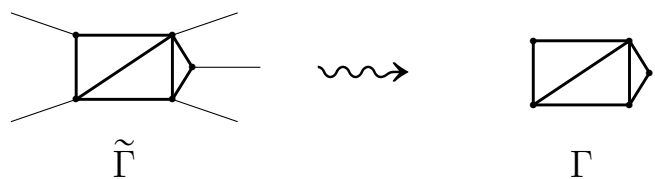

\paragraph{Formal setup: momenta.}\label{par:formal-momenta}
Let $V$ be a  
 $\CC$-vector space  of dimension $D\geq 3$
 with a non-degenerate symmetric $\mathbb{C}$-bilinear form $(-,-)$. Although this form identifies $V$ with $V^{*}$, we will think of $V$ as the space of momenta in (i.e., the dual vector space to) the complex Minkowski space $V^{*}$. For $k \in V$ we write $k^{2}=(k, k)$. 
 
  The space $V^L=\bigoplus_{l\in L} V$ consists of all assignments $p=(p_l)_{l\in L}$
  of momenta
 to the legs of $\wt\Gamma$.

 \vskip .2cm
 
 Let  $C_\bullet(\Gamma, V)$ and $C_\bullet (\wt\Gamma, V)$
  be  the chain complexes of $\Gamma$  and $\wt\Gamma$ with
coefficients in $V$.  
Each of them reduces to a single map:
\be\label{eq:chains-Gamma}
\xymatrix{
\bigoplus_{e\in\Ed(\Gamma)} \orr(e)\otimes V = C_1(\Gamma, V)   \ar^\del[rr] && 
C_0(\Gamma, V) = \bigoplus_{v\in \Vert(\Gamma)} V,}
\\
\xymatrix{
\bigoplus_{e\in\Ed(\wt \Gamma)} \orr(e)\otimes V = C_1(\wt \Gamma, V)   \ar^{\wt\del}
[rr] && 
C_0(\wt \Gamma, V) = \bigoplus_{v\in \Vert(\Gamma)} V = C_0(\Gamma, V). 
}
\ee
Here $\orr(e)$ is the $1$-dimensional {\em orientation space} of $e$. That is, 
a choice of
orientation of $e$ identifies $\orr(e) \to \CC$ and a different choice changes the
sign of the identification. As our graphs are not allowed to have loops,
any choice of a vertex $v$  (out of one or two) of $e$ defines an orientation of $e$
``towards $v$''. So an element $\sigma$ of, say, $C_1(\wt\Gamma, V)$
can be seen as a collection of vectors $\sigma(v,e)\in V$ associated to all
incident pairs $(v,e)$ of of a vertex and edge, 
so that whenever $e$ is incident to another vertex $v'$, we have 
$\sigma (v', e) = - \sigma(v,e)$.

\vskip .2cm


Note that $\wt\Gamma$ is not compact (legs are open intervals) and
$C_\bullet(\wt\Gamma, V)$ calculates $H_\bullet^{\on{BM}}(\Gamma, V)$,
the {\em Borel-Moore}, or {\em  locally finite} homology, dual to 
the  cohomology with compact support.
This homology is not homotopy invariant, depending on the number of legs. 
It is useful to introduce the following terms. 

\begin{defi}\label{def-scat-corolla}
(a) Let $I$ be a finite set. A {\em scattering corolla} labelled by $I$ is a
collection $p=(p_i)_{i\in I}$ of vectors from $V$ such that $\sum p_i=0$.
The space of such collections will be denoted $V^{I,0}\subset V^I$. 

\vskip .2cm

(b) Let $\wt\Gamma$ be a graph as before. A {\em scattering diagram}
labelled by $\wt\Gamma$ is a $1$-cycle $\sigma$  in $C_1(\wt\Gamma, V)$,
i.e., an element in $\Ker(\wt\del)$. 
\end{defi}

Thus a scattering corolla represents a collection of momenta satisfying the
conservation law, while a scattering diagram represents
an assignment of momenta to edges of $\wt\Gamma$ such that
the conservation law holds at each vertex. This means that
a scattering diagram  $\sigma$ gives rise to several scattering corollas:
\begin{itemize}
\item The {\em external} scattering corolla formed by the  momenta
$p_l = \sigma(\eps(l), l)$ associated to legs $l\in L=\Leg(\wt\Gamma)$. 

\item The {\em local} scattering corolla $\sigma_v$  at each vertex $v$. 
It is  formed by the momenta
$k_e = \sigma(v,e)$ associated to all the edges  $e$ (including legs) incident to $v$
and oriented towards $v$. 
\end{itemize}

\noindent We can also say that a scattering corolla is a scattering diagram
labeled by the {\em corolla graph} $\Cor(I)$ with one vertex and $I$ legs.
We have the natural projection
\be
\Pr: \wt\Gamma \lra \Cor(L)
\ee
which contracts the legless part $\Gamma$ into the unique vertex.
 Cf.  Fig. \ref{fig:S-and-graphs},
where the  left part (with the blob labelled $S$) represents $\Cor(L)$. 
The external scattering corolla of a scattering diagram $\sigma$
can be written as $p=\Pr_*(\sigma)$, the direct image of the cycle
$\sigma$ under $\Pr$.   We will say that the scattering diagram $\sigma$
{\em refines} its external  scattering corolla $p$.

\vskip .2cm

More explicitly, recall that 
 $\Ed(\wt\Gamma)  = L\sqcup \Ed(\Gamma)$. Note also that each leg $l$
 of $\wt\Gamma$ has a canonical orientation (towards the unique  vertex $\eps(l)$). 
Further, 
the map $\eps: L\to\Vert(\Gamma)$
defines a linear map 
\[
\eps_*: V^{L,0} \to C_0(\Gamma, L),
\quad 
(\eps_* (p))_v \,=\, - \sum_{\eps(l)=v} p_l
\]
(the minus sign added for convenience)
and we can identify $C_\bullet(\wt\Gamma, V)$ as 
\be
\xymatrix{
 V^L \oplus  C_1(\Gamma, V)   \ar^{\hskip 0.5cm \wt\del\,\, =\,\,   -\eps_* + \del}
[rr] && 
  C_0(\Gamma, V). 
}
\ee
So a scattering diagram $\sigma$ for $\wt\Gamma$ is the same as a pair
$(p, k)\in V^{L,0} \times C_1(\Gamma, L)$ such that $\eps_*(p) = \del(k)$.
In other words,  for a  scattering corolla $p\in V^{L,0}$ the subspace $\del^{-1}\eps_* (p)$
 consists of  scattering diagrams labelled by $\wt\Gamma$ 
 for which $p$ is the external corolla.

 \vskip .2cm
 
 Denote by $b= b_1(\Gamma)$ the first Betti number (the number of loops) of $\Gamma$.
 As $\Gamma$ is connected,  the image $\Im(\del)$  consists of $q=(q_v)\in C_0(\Gamma, V)$
 such that $\sum q_v=0$. In particular, the image of $\eps_*$ is contained in $\Im(\del)$.
 As for the kernel $\Ker(\del)$, it has has dimension $bd$, being identified with
 $H_1(\Gamma, V) = H_1(\Gamma, \CC) \otimes V$. Thus $\del^{-1}\eps_*(p)$ is always a complex affine space of
 dimension $bD$.

 \paragraph{Formal setup: integrals.} 

Each edge
 $e \in \operatorname{Ed}(\Gamma)$  represents an``intermediate particle"  of mass $m_{e}$ which has  the "propagator" (Green's function in the momentum
 representation) which is a function on $V$.  Since we consider only scalar particles,  the propagator has the form
\[
k \mapsto \frac{1}{A_{e}(k)} , \quad A_{e}(k)= m_{e}^{2}-k^2.
\]
Note that $A_e$ is an even function on $V$ and  each torsor $\orr(e)$
is trivialized up to sign $\pm 1$. So we can consider $A_e$ as a function on 
$\orr(e)\otimes V$ and, by pullback, on $C_1(\Gamma, V)$. 

\vskip .2cm

Next, the interaction terms in the Lagrangian give, at each vertex $v$, a polynomial 
$B_{v}(p,k)$ depending on the momenta $p_l$,  $k_{e}$ of the particles corresponding to legs and edges  incident to $v$. This produces  an algebraic
 function (rational, if $D$ is even)
\be\label{eq:integrand}
\begin{gathered}
F: C_1(\wt\Gamma, V) = V^{L,0}\times C_{1}(\Gamma, V) \longrightarrow \mathbb{C}, 
\\
\bigl(p = (p_l)_{l\in L}, 
 k=\left(k_{e}\right)_{e \in \mathrm{Ed}(\Gamma)}\bigr)
 \,\,\mapsto \,\,  F(p,k) = \frac{B(p,k)}{\prod_{e} A_{e}\left(k_{e}\right)}, \\
   B(p,k)=\prod_{v} B_{v}(p,k) .
  \end{gathered} 
\ee
 So $B(p,k)$ is a polynomial. The {\em Feynman integral} corresponding to 
 $\wt\Gamma$ is  the function $I=I_{\wt\Gamma}$ on $V^{L,0}$ defined by
\be\label{eq:integral} 
p \mapsto I(p) =  I_{\wt\Gamma}(p)\,=\, \int\limits_{k\,\in\, \gamma(p)\,\subset\,  
\partial^{-1}\eps_*(p)} F(p,k) d k, \quad p\in V^{L,0}, 
\ee
where the integration is taken over some real cycle $\gamma = \gamma(p)$ of dimension $bD$
inside $\partial^{-1}\eps_*(p)$
avoiding singularities of the integrand and   varying with $p$ according to the Gauss-Manin connection.  
That is, in the diagram
\be\label{eq:diagram-push-pull}
\xymatrix{
& V^{L,0}\times C_1(\Gamma, V)  = C_1(\wt\Gamma, V)
\ar[d]^{\Id\times\del}
\\
V^{L,0} \ar[r]_{(\Id, \eps_*)\hskip 1cm} & V^{L,0}\times C_0(\Gamma, V)
}
\ee
we fist take the direct image  of $F(p,k)$
under the vertical map (i.e., integrate it along  the fibers)  and then pull the result
back to $V^{L,0}$  via the horizontal map. 

\vskip .2cm

As there is no canonical choice of $\gamma(p)$, the function $I(p)$
is, in general, multivalued. 
We want to find out the locus of $p$ where it can be singular (go to infiity or ramify). 


\section{Landau conditions: complex domain}\label{sec:landau-complex}

\paragraph{ Landau patterns. }
We have the hypersurface $Y \subset C_{1}(\Gamma, V)$ given by $\prod_{e} A_{e}\left(q_{e}\right)=0$. So $Y=\bigcup_e Y_e$ is reducible with the component $Y_e$ given by
$A_e(q_e)=0$ lifted from the $e$th factor.
 
 As the integration cycle $\gamma(p)$ in \eqref{eq:integral} lies in
$\partial^{-1}\eps_*(p) \- Y$, singularities of $I(p)$ happen when:
\begin{itemize}
\item The number of independent integration cycles changes, i.e.:

\item The topology of $\partial^{-1}\eps_*(p)\- Y$ changes, i.e.:

\item The topology of $Y\cap \partial^{-1}\eps_*(p)$ changes (compared to generic $p$). 
\end{itemize}

\noindent Note that this does not depend on the numerator $B(q)$, i.e.,
on the details of the interactions, but only on the propagator. This is
the first observation  of Landau \cite{landau1}. 

\vskip .2cm

A precise modern formulation of the above conditions would be in terms
if microlocal geometry \cite{KS}. That is, we should form the
Lagrangian variety  $\Lambda \subset T^* C_1(\Gamma, V)$ as the union of the closures
of the conormal bundles to the smooth loci of the $Y_e$
 and then
apply to it the operations of Lagrangian direct and inverse images
in the diagram \eqref {eq:diagram-push-pull},
as defined in \cite{KS}. This was analyzed in detail in \cite{sato}. 
In this paper we would like to focus on the most immediate 
type of scenario highlighted in \cite{landau1}. 

\vskip .2cm

Denote $Y_p = Y\cap  \partial^{-1}\eps_*(p)$ (the singular locus to be avoided
in evaluating $I(p)$). Then $Y_p = \bigcup_e Y_{p,e}$ where $Y_{p,e} = Y_e\cap
\del^{-1}\eps_*(p)$ is a quadric hypersurface. 

\begin {defi}
A {\em Landau pattern} consists of:
\begin{itemize}
\item[(1)] A subset $I\subset \Ed(\Gamma)$;

\item[(2)] A collection of momenta $q=(q_e\in V), e\in I$ such that
$q_e\in Y_{e,p}$ (``on-shell'' condition) and the differentials $d_{q_e}(A_e)$
are linearly dependent as functions on $T_q(\del^{-1}\eps_*(p))$
but any proper subset of these differentials is linearly independent. 
\end{itemize}
Here $q$ is considered as a point of $C_1(\Gamma, V)$ by
putting $q_e=0$ for $e\notin I$. 
\end{defi}

In other words, in a Landau pattern we have
\be\label{eq:land-dependence}
\sum_{e\in I} \alpha_e \,\cdot d_{q_e} A_e|_{T_q(\del^{-1} \eps(p))} = 0
\ee
for some collection $\alpha=(\alpha_e)_{e\in I}$ of  nonzero complex numbers. 
These numbers are called  the {\em Landau multipliers}. 

\vskip .2cm

Note that the requirement that all $\alpha_e\neq 0$ (i.e., no proper
subset of the $d_{q_e} A_e, e\in I$ is linearly dependent) is no
restriction of generality since we can always restrict $I$ to include only
$e$ corresponding to a minimal dependent subset 
(i.e, those with $\alpha_e\neq 0$).

\vskip .2cm

Clearly, a Landau pattern defines a singular point $p$ for $I(p)$
(unless this pattern persists of all $p\in V^L$), as the differentials of the equations
of the hypersurfaces $Y_{e,p}$ become linearly dependent somewhere
on the intersection of these hypersurfaces. 
We will refer to singular points thus obtained as {\em Landau singularities},
understanding the term in a stricter sense than, say, \cite{pham}.

\begin{rem}
Note that in   Feynman's  diagram calculus the intermediate particles are not, a priori,
subject to  the on-shell condition $q_e^2=m_e^2$. However, if
we  are interested in 
  singularities of $I(p)$ (of the above  most immediate type), only on-shell
intermediate particles contribute to them. 
This was the second observation of Landau \cite{landau1}. 
\end{rem}

\paragraph{The Landau conditions: spatial  realizations of $\Gamma_I$.} \label{par:spatial-real}
To analyze the condition \eqref{eq:land-dependence} for a Landau pattern,
we write the differential of $A_e$ as
 \[
d_{q_{e}}\left(A_{e}\right)=d_{q_{e}}\left(q \mapsto m_{e}^{2}-q^2\right)=-2(q_{e}, d q).
\]
In other words, the value of $d_{q_{e}}\left(A_{e}\right)$ on a tangent vector $\dot{q}=\left(\dot{q}_{e}\right) \in T_{q}( \partial^{-1}\eps_*(p))$ is equal to
 $-2 (q_{e}, \dot{q}_{e})$. 
 Further, the tangent space to the affine space 
$\partial^{-1}\eps_*(p)$ is the vector space $\partial^{-1}(0)=Z_{1}(\Gamma, V)$ 
 of {\em 1-cycles}  of $\Gamma$ with coefficients in $V$. 

\vskip .2cm

Consider now the {\em cochain complex} of $\Gamma$ with coefficients
in $V$:
\be
 \xymatrix{
 C^0(\Gamma, V) = \bigoplus_{v\in \Vert(\Gamma)} V  \ar^{\delta\hskip 1cm}[rr] && 
 \bigoplus_{e\in\Ed(\Gamma)} \orr(e)\otimes V = C^1(\Gamma, V). 
}
\ee
Its terms are formally identical to those of the chain
complex \eqref{eq:chains-Gamma}
 but we use the
Minkowski scalar product $\langle - , - \rangle$ to identify
$C^i(\Gamma, V)$ with the dual of $C_i(\Gamma, V)$.
Then the coboundary 
map $\delta$ becomes dual to $\del$. 
\vskip .2cm

 Let us consider the tuple
$  (\alpha_e q_e)_{e\in\Ed(\Gamma)}$ as an element of $C^1(\Gamma, V)$
(with components at $e\notin I$ being $0$). The condition
 \eqref{eq:land-dependence} means that it is orthogonal
to the space of cycles $Z_1(\Gamma, V)$. Because of the duality
between $\del$ and $\delta$, this   means that $  (\alpha_e q_e)$
is 
a {\em coboundary}, i.e., lies in the image of $\delta$. 
 In other words, there is a map
$\phi: \Vert(\Gamma) \to V$ such that for any edge $e$ with vertices $v,w$
(oriented from $v$ to $w$) the component $\alpha_e q_e$
equals $\phi(w)-\phi(v)$, in particular, $\phi(v)=\phi(w)$ if $e\notin I$.  
In other words, edges of $\Gamma$ become represented by
straight intervals in $V$. 
This leads to the following reformulation ({\em Landau conditions}). 

\vskip .2cm

Let $\Gamma_I$ be the quotient graph of $\Gamma$ obtained by
contracting all the edges not in $I$, so $\Ed(\Gamma_I)=\Ed(\Gamma) \- I$ and
$\Vert(\Gamma_I)$ is a quotient set of $\Vert(I)$.

\begin{prop}\label{prop:landau-reform}
A collection $q = (q_e)_{e\in I}$ represents a Landau pattern
if and only if there exist nonzero $\alpha_e\in \CC$, $e\in I$ 
 and a piecewise-affine map $\bphi: \Gamma_I \to V$ that $\alpha_eq_e$
 is the value of $\phi$ on the edge $e$. That is, if $e$ is oriented
 from a vertex $v$ to the other vertex $w$, then 
 $\alpha_eq_e = \bphi(w)-\bphi(v)$. \qed
\end{prop}

We will  call such $\bphi$  a {\em spatial realizations} of $\Gamma_I$. 
It is uniquely defined by its values on vertices, i.e., by a $0$-cochain
$\phi\in C^0(\Gamma, V)$ and we will use the two interchangeably.
Note that a spatial realization is injective on all edges but does not
 have to be  an embedding of $\Gamma_I$ into $V$.  For example,
   two vertices
not connected by an edge may be sent to the same point of $V$
and also, the images of different edges may intersect somewhere  in their 
intermediate
points. 

\paragraph{Principal Landau singularities.} Landau singularities
corresponding to Landau patterns with $I=\emptyset$
(so that $\alpha_e\neq 0$ for all $e\in\Ed(\Gamma)$) will
be called {\em principal Landau singularities}. They correspond to
spatial realizations of $\Gamma$ itself. Landau singularities corresponding to 
nonempty $I$ are pullbacks of principal Landau singularities for
quotient graphs of $\Gamma$. We note that vertices of such quotient
graphs are agglomerations of vertices of $\Gamma$ and
so may not correspond to actual monomials in our original Lagrangian. 
Nevertheless, as far as singularities are concerned, they are pullbacks
and so in the sequel we will restrict our attention to principal
Landau singularities only. 

\paragraph{Producing Landau singularities geometrically.} 
\label{par:landau-geom} 
Reversing the logic, we frame the above analysis as a geometric
way of constructing (principal) Landau singularities. It is remindful
of Hodge theory on $\Gamma$ in that, identifying chains and cochains,
we consider both differentials simultaneously:
\be
 \xymatrix{
 C^0(\Gamma, V)    \ar@<.4ex>[rr]^{\delta=\del^*} && 
  \ar@<.4ex>[ll]^{\del}
  C^1(\Gamma, V). 
}
\ee
 It consists of the following steps, some of which
 restrict the previous ones:
 \begin{itemize}
 \item[(1)] Start with a spatial realization $\bphi$ of $\Gamma$ 
 such that edges $e$ with $m_e=0$ are sent to lightlike (null) intervals. 
 We can identify $\bphi$ with a $0$-cochain $\phi\in C^0(\Gamma, V)$. 
 
 \item[(2)] Form the {\em coboundary} $\delta\phi$, so that $(\delta\phi)_e$
 is the vector corresponding to the image of $e$ under $\bphi$. 
 
 \item[(3)] Stretch each $(\delta\phi)_e$ by a nonzero $\alpha_e\in\CC$
 so that $(\alpha_e^{-1} (\delta\phi)_e)^2 = m_e^2$. This allows
 exactly two possibilities for $m_e\neq 0$ and any $\alpha_e$
 is allowed for $m_e=0$. This gives a new cochain $\alpha^{-1} \cdot\delta\phi$. 
 
 \item[(4)] Now form the {\em boundary}  
 $\del (\alpha^{-1} \cdot \delta\phi)\in C^0(\Gamma, V)$ which we view as
 an assignment of momenta to vertices of $\Gamma$. 
 The choice of $\phi$ above should be such that $\del (\alpha^{-1}\cdot \delta\phi)$
 lies in the image of $\eps_*: V^L\to C^0(\Gamma, V)$. In this
 case any $p\in V^L$ with $\eps_*(p) = \del (\alpha^{-1} \cdot \delta\phi)$
 represents a Landau singularity. In other words, we can distribute
 $\del (\alpha^{-1} \cdot \delta\phi)$ among the incoming legs arbitrarily. 
 \end{itemize}
 
 \begin {rem}
 In particular, we see that a Landau singularity $p\in V^{L,0}$,  
 a scattering corolla (in the sense of Definition \ref {def-scat-corolla}), 
 is always realized as the external corolla of 
 a scattering diagram, namely  of  $(p, q)$ where $q = \alpha^{-1}\cdot \delta\phi$. 
 \end{rem}

 \noindent Continuing the analogy with Hodge theory, we see
 that  the key role in Landau's construction is played by the operator
 \be
 \DD \, = \del \circ \alpha^{-1}\circ \delta \,=\,\delta^* \circ\alpha^{-1} \circ\delta:
 \,\,\,  C^0(\Gamma, V) \lra 
 C^0(\Gamma, V), 
 \ee
 where $\alpha^{-1}$ is the diagonal operator multiplying the summand $V$ associated
 to $e$ by $\alpha_e^{-1}$.
 The operator $\DD$ can be seen as the {\em discrete Laplacian}  associated to $\Gamma$.
 The datum $\alpha = (\alpha_e)_{e\in \Ed(\Gamma)}$ plays the role of the
 metric on the graph.
 
\begin{figure}[h]
 \centering
 \begin{tikzpicture}[scale=0.2]
 
 \node at (0,0){\small$\bullet$}; 
  \node at (14,0){\small$\bullet$}; 
 \node at (0,14){\small$\bullet$};  
  \node at (8,12){\small$\bullet$}; 
   \node at (14,7){\small$\bullet$}; 
   
   \draw[line width=0.7] (0,0) -- (0,14) -- (8, 12) -- (14,7) -- (14,0) -- (0,0); 
      \draw[line width=0.7]  (0,14) -- (14,0); 
         \draw[line width=0.7] (0,0) -- (8, 12); 
         
  \draw (0,0) -- (-4, -8) ;
  \draw (0, 14) -- (-2, 18); 
  \draw (8, 12) -- (24, 21);  
  \draw (14,7) -- (21, 21); 
  \draw (14,0) -- (20, -3);

  \draw[->, line width=1] (0,0) -- (0,2);  
    \draw[->, line width=1](0,0) -- (2,0); 
      \draw[->, line width=1] (0,0) -- (-2, -4);   
   \draw[->, line width=1] ((0,0) -- (2,3); 

 \node at (-5.5, -7){\small$\wt\bphi(l)$};
\node at (-3.5, -2.5){$p_l$}; 

\node at (-3.5, 0){\small$\wt\bphi(\eps(l))$}; 
  
  \node at (24, 7){$\wt\bphi(\wt\Gamma)$}; 
   
 \end{tikzpicture}
 \caption{A Landau graph:  the images of edges (including legs) can intersect.} 
 \label{fig:landau-graph}
 \end{figure}
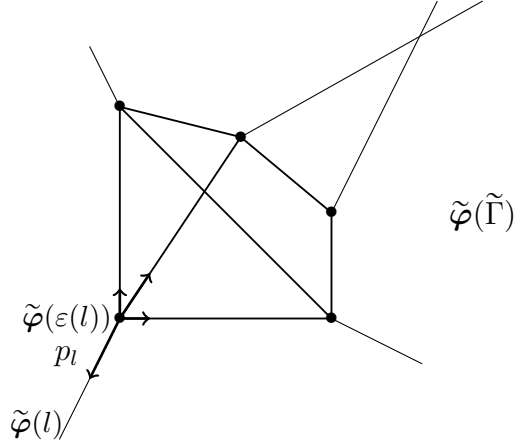

 \paragraph{Landau graphs.} It is convenient to extend the spatial realization
 of $\Gamma$ to that of $\wt\Gamma$ by adding the images of the legs
  encoding the incoming momenta. 
  That is, we map the leg $l\in L$  with vertex $v=\eps(l)$ to the real half-line
 $\phi(v) + \RR_+\cdot p_l$. 
 This defines a piecewise-affine map
 $\wt\bphi: \wt\Gamma\to V$ extending $\bphi$.. The map $\wt\phi$ is
  an {\em immersion}, i.e., its restriction  to each sufficiently small part
  of $\wt\Gamma$ is an embedding, but is need not be an embedding globally,
  see Fig. \ref {fig:landau-graph}. For future use, we introduce the following

\begin{defi}\label{def:landau-graph}
A {\em Landau graph}  is a datum   of a colored graph $\wt\Gamma$ as above
and a piecewise-affine immersion $\wt\bphi: \wt\Gamma\to V$ such that
there exists a scattering diagram 
$\sigma =
\bigl( (p_l)_{l\in L}, (q_e)_{e\in\Ed(\Gamma)}\bigr)$
with the properties:
\begin{itemize}
\item[(1)] Each $p_l\in V$ lies in the real positive direction of the half-line
$\wt\bphi(l)$.

\item[(2)] Each $q_e\in V\otimes\orr(e)$ lies in the complex direction of
the edge $\wt\bphi(e)$. 
\end{itemize}
Such $\sigma$ will be called {\em collinear} with $(\wt\Gamma, \wt\bphi)$. 
\end{defi}

 More precisely, (1) means that $p_l \in \RR_{>0}\cdot (a-\bphi(\eps(l))$ for any
 $a$ lying on the half-line $\wt\bphi(l)$, other than the vertex $\wt\bphi(\eps(l))$.
 Similarly,  (2) means that 
 \[
 q_e\in \bigl(\CC^*\cdot (\wt\bphi(w) -\wt\bphi(v)\bigr)
 \otimes \orr(e), \quad
 \text{where} \quad  e=[v,w].
 \]  
 Fig. \ref{fig:landau-graph} depicts a Landau graph and one local corolla of a
 collinear scattering diagram.  Proposition
 \ref {prop:landau-reform} can be reformulated geometrically as follows.
 
 \begin{prop}\label{prop:Lsing=Lgrap}
 A collection $p=(p_l)_{l\in L}$ represents a (principal) Landau singularity
 for $\wt\Gamma$
 if and only if it  comes from a Landau graph $(\wt\Gamma, \wt\phi)$, i.e., $p$ is 
 the external corolla of a scattering diagram $\sigma$ 
 collinear with $(\wt\Gamma, \wt\phi)$. \qed
 
 \end{prop}
 
 
 \section{Landau conditions in context}\label{sec:landau-context}
 
 Here we recall some other classical constructions related to
 the Landau conditions. 
 
 \paragraph{Relation to the Feynman parametrization and the Cayley trick.} 
 Landau's analysis was based on the so-called {\em Feynman parametrization} which
 reduces integration of a product $f_1(x)^{-s_1}\cdots f_n(x)^{-s_n}$ of powers
 of several functions to integration of a power of a single function
 \[
 F(\alpha, x) \,=\,\sum_{i=1}^n \alpha_i f_i(x)
 \]
 involving additional variables $\alpha_i$, times some powers of these variables. 
 This  is based on the classical
 identity holding for any nonzero numbers $f_i\in\CC^*$ and any $s_i\in \CC$ with $\Re(s_i) >0$:
 \[
 {1\over f_1^{s_1} \cdots f_n^{s_n}} = 
 {\Gamma(s_1+\cdots  s_n)\over 
 \Gamma(s_1) \cdots \Gamma(s_n)}\,\, 
 \int_{\alpha \in \Delta^{n-1}} 
 {
 \alpha_1^{s_1-1} \cdots \alpha_n^{s_n-1} \over
 (\alpha_1 f_1 + \cdots + \alpha_n f_n)^{s_1+\cdots + s_n} 
 } \, d\alpha.  
\]
Here $\Delta^{n-1}$ is the simplex
\[
 \Delta^{n-1} = \bigl\{ \alpha = (\alpha_1, \cdots, \alpha_n) \in \RR^n\,\bigl| \, 
 \alpha_i \geq 0, \,\,\sum\alpha_i = 1 \bigr\}  \quad\text{and} \quad  d\alpha = d\alpha_1 \cdots d\alpha_{n-1}. 
 \]
 In our particular case  where  all  $s_i= 1$ and $n=|\Ed(\Gamma)|$ this gives
 \be\label{eq:I(p)-Feynman}
 I(p) \,=\,
 (n-1)! 
 \int\limits_{(k, \alpha)\, \in \, \gamma(p) \times\Delta^{|\Ed(E)|-1} }
 {
 B(p,k)
  \over F(\alpha, k)^{|\Ed(\Gamma)|} 
 } dk \, d\alpha,
 \ee
 where
 \be\label{eq:F-alpha-k}
 F(\alpha, k) \,=\,\sum_{e\in\Ed(\Gamma)} \alpha_e \, A_e(k_e)
 \ee
 is an inhomogeneous  quadratic function on $C_1(\Gamma, V)$ depending linearly on the $\alpha_e$,
 which we consider in  its  restriction to $\del^{-1}\eps_*(p)$.
 
 \vskip .2cm
 
 In algebraic geometry these considerations take form of the following general principle
 which can be traced back to Cayley and which was called in \cite{GKZ-Euler, GKZ} the
 {\em Cayley trick}:
 
 \vskip .2cm
 
 {\em 
 Degeneration (non-generic behavior) of a configuration of hypersurfaces 
 $Y_i = \{f_i(x)=0\}$, $i=1,\cdots, n$, 
 inside some smooth algebraic variety $X$ is equivalent to degeneration of a single hypersurface
 $\{\sum \alpha_i f_i(x)=0\}$ inside $X\times \CC^n$ considered together with the coordinate hyperlplanes
 $\alpha_i=0$.
 }
 
 \paragraph{The Symanzik polynomials.} For completeness let us
 explain the conceptual meaning of the classical Symanzik polynomials
  \cite{bloch, mizera, nakanishi, smirnov} in our notation with $C_\bullet(\Gamma, V)$. 
 Let
 \[
 F(x) = (\ba x, x) + (\bb, x) + c,
 \,\,\, F^{(2)} (x) = (\ba x, x), \quad x=(x_1,\cdots, x_m)^t
 \]
 be an inhomogeneous quadratic function in $m$ variables $x_1,
 \cdots, x_m$ and its homogeneous
  quadratic part. So $\ba = \|a_{ij}\|$ is a symmetric $m\times m$ matrix,
 $\bb = (b_1,\cdots, b_m)$ is a covector and $c$ is a number. 
 We consider these coefficients as indeterminate variables. 
 
 \vskip .2cm
 
 The {\em discriminant} of $F$ can be defined, up to
 a constant factor, as an  irreducible polynomial
 $\Disc(F) = \Disc(a_{ij}, b_i, c)$ in the coefficients which
 vanishes when the hypersurface $F=0$ is singular. In other words,
 $\Disc(F)$ is the $A$-discriminant \cite{GKZ} where $A$ is the
 set of monomials in the $x_i$ of degree $\leq 2$.  To find
 $\Disc(F)$ explicitly, we homogenize $F$, converting it to a quadratic
 form in $m+1$ variables
 \be\label{eq:quad-homog}
 \wt F(y_0, \cdots, y_m) \, =\, y_0^2 \, F(y_1/y_0, \cdots, y_m/y_0)
 \ee
 and take the determinant of the corresponding
 $(m+1)\times (m+1)$ matrix:
 \[
 \Disc(F) = \det \begin{pmatrix} \ba & {1\over 2} \bb^t
 \\
 {1\over 2} \bb & c\end{pmatrix} \, = \, (c-{1\over 4} \bb\,  \ba^{-1}\,
  \bb^t)\cdot \det(\ba).
 \]
 The second equality is an instance of a general identity for
 the determinant of any block $2\times 2$ matrix, cf.
 \cite{gelfand-retakh}.  For $m=1$ when $\ba=a, \bb=b, c$ are
 all numbers, this gives $ac-{1\over 4} b^2$ which is $-{1\over 4}$ times  
 the standard form $b^2-4ac$. 
 
 \vskip .2cm
 
 Let also 
 \[
 \Disc(F^{(2)}) = \det(\ba)
 \]
 be the discriminant of the quadratic form $F^{(2)}$, i.e., the $A$-discriminant
 with $A$ being the set of monomials in the $x_i$ of degree exactly $2$. 
 The expression
 \[
 \Crv(F)\, =\,  {\Disc(F)\over\Disc(F^{(2)}) } \,= \,
 c-{1\over 4} \bb\,  \ba^{-1}\, \bb^t
 \]
 is the {\em critical value} of $F$, i.e., its value at the (unique, generically)
 critical point $x= -{1\over 2}   \ba^{-1}\bb^t$. 
 Both $\Disc(F)$ and $\Disc(F^{(2)})$ are unchanged under the translation of
 variables in $F$. 
 
  \vskip .2cm
 
 We apply this to the family of quadratic functions $F(\alpha, k)$ of 
 \eqref{eq:F-alpha-k} on the variable affine space $\del^{-1}\eps_*(p)$,
 $p\in V^{L, 0}$. As each $\del^{-1}\eps_*(p)$ is identified with
 $\del^{-1}(0)$ up to a translation, the discriminant of
 $F(\alpha, k)|_{\del^{-1}\eps_*(p)}$ makes sense as a polynomial in
 $\alpha$ and $p$. In this notation the first and second Symanzik
 polynomials are
 \be
 \begin{gathered}
 \Uc(\alpha) \,=\,\Disc\bigl(F(\alpha, k)^{(2)}|_{\del^{-1}\eps_*(p)}\bigr) =
 \Disc\bigl(F(\alpha, k)^{(2)}|_{\del^{-1}(0)}\bigr),
 \\
 \Vc(\alpha, p) \,=\, \Disc\bigl( F(\alpha, k)|_{\del^{-1}\eps_*(p)}\bigr).   
 \end{gathered}
 \ee
 The ratio
 \be
 {\Vc(\alpha, p)\over \Uc(\alpha)} \,=\, \Crv\bigl(F(\alpha, k)|_{\del^{-1}
 \eps_*(p)}\bigr),
 \ee
 i.e.,  the critical value of $F(\alpha, k)$ on $\del^{-1}
 \eps_*(p)$, is the quantity denoted by $\phi$ in
 \cite[\S2]{landau1}.  
 
 \paragraph{Analogy with hyperdeterminants.} Denoting by $(x_1, \cdots, x_m)$,
 $m=bD$, 
 some linear coordinates on $\del^{-1}\eps_*(p)$, we homogenize
 $F(\alpha, k)$  to a function $\wt F(\alpha, y)$ with $y=(y_0, \cdots, y_m)$
 as in \eqref{eq:quad-homog}. This gives a function 
 \[
 \wt F(\alpha, y) \,=\,\sum_{e\in\Ed(\Gamma)} \sum_{j,l=0}^m  c_{ejl}(p) \alpha_e y_j y_l
 \]
 linear homogeneous
  in $\alpha$ and
 quadratic homogeneous in $y$ so that  $p$ is a Landau singularity when
 this function has a critical point $(\alpha^{(0)}, y^{(0)})$ with all $\alpha^{(0)}_i\neq 0$
 and $y^{(0)}\neq 0$. The problem of finding the condition for  such functions to
 have a 
 nontrivial critical
 point is a version of the problem of the  {\em hyperdeterminant} \cite{GKZ} that
 concerns a trilinear function
 \[
 \begin{gathered}
 T(x,y,z)  = T_B(x,y,z)  = \sum_{i=1}^p \sum_{j=1}^q \sum_{l=1}^r \,  b_{ijl} \,
 x_i, y_j z_l,
 \\
 x=(x_1, \cdots x_p), \,\, y= (y_1, \cdots, y_q), \,\, z= (z_1, \cdots, z_r)
 \end{gathered}
 \]
 associated to a $3$-dimensional matrix $B=\|b_{ijl}\|$ of size $p\times q\times r$.
 The locus of $B$ such that  $T_B$ has a nontrivial critical point is called the
 {\em hyperdeterminantal variety}  and the hyperdeterminant proper is the equation
 $\Det(B)$ of this variety in  the case when it has codimension $1$. 
 
 \vskip .2cm
 
 In our case we have a version of this situation when
   $p=|\Ed(\Gamma)|$,  $q=r=bD+1$, the matrix $C=\|c_{ejl}\|$
 is symmetric in $j,l$ but  we consider $T_C(\alpha, y,y)$ instead of
 $T_C(\alpha, y,z)$. These two settings: a trilinear  vs. a linear-quadratic
  form,  are not exactly identical but closely related
  and some relations
 between them were discussed in \cite[\S8]{ottaviani}.
 
  \vskip .2cm
  
  From this point of view, the concept of the Landau discriminant
  introduced in \cite {mizera}  can be seen as an analog of the 
  Schl\"afli method for computing hyperdeterminants \cite{GKZ}.


\section{Real domain:   generalities}\label{sec:real-gen}

\paragraph{The real setup: positivity of the $\alpha_e$.}\label{par:alpha>0} 
 We now restrict to the real Minkowski
space $V_\RR= \RR^D\subset V$ with the real  non-degenerate scalar product 
$(  -,- )$. We consider two possibilities:
\begin{itemize}
\item[(A)] $(  -,- )$ has Lorenzian signature  $(+ -\cdots -)$. 
This is the standard physical assumption.

\item[(B)] $(  -,- )$ is  Euclidean, i.e., of signature $(+\cdots +)$. 
This is also a common assumption, corresponding to spatial directions
being ``imaginary''. 

\end{itemize}
 The chain and cochain differentials
being defined over $\RR$ (in fact, over $\ZZ)$, we have the real 
operators
\be
 \xymatrix{
 \bigoplus_{v\in \Vert(\Gamma)} V_\RR \,\,= \,\,  C^0(\Gamma, V_\RR)    \ar@<.4ex>[rr]^{\delta=\del^*} && 
  \ar@<.4ex>[ll]^{\del}
  C^1(\Gamma, V_\RR) \,\,= \,\, \bigoplus_{e\in\Ed(\Gamma)} V_\RR \otimes
  \orr(e)
  \\
  V_\RR^{L,0} \,\,= \,\,  V^{L, 0}\cap   \bigoplus_{l\in L} V_\RR.   \ar[u]^{\eps_*} &
}
\ee
We are interested in Landau singularities of $I(p)$ lying in $V_\RR^{L,0}$,
the set of collection of real momenta summing up to $0$. As there
are many possible branches of $I(p)$ (they corrrespond to
various cycles $\gamma(p)$), we focus on some particular ones
which are immediately relevant to the real setup and look at
 their singularities only. 

\vskip .2cm

First, we note that for $p\in V_\RR^{L,0}$ real the affine space
$\del^{-1}\eps_*(p)$ is also real and the (inhomogeneous)
quadratic functions $A_e(q_e)$
on this space are real as well. So it is natural (and we will do so)
 to focus on
cycles $\gamma(p)$ (in the complexification of  $\del^{-1}\eps_*(p)$)
 which are as close as possible to the real locus
 (e.g., avoiding real singularities by some kind of $i\eps$ prescription)
  and on the corresponding
 branches of $I(p)$. A Landau singularity in such a branch would correspond
 to a Landau pattern with all the data   being real, i.e., $q_e\in V_\RR$
 and $\alpha_e\in \RR^*$.  This allows $2^{|\Ed(\Gamma)|}$
 possibilities for the signs of the $\alpha_e$. However, we make
 the following more restrrictive
 
 \begin{defi}\label{def:real-landau}
(a)  A {\em real Landau pattern} is a Landau pattern with all $q_e\in V_\RR$ and
 all $\alpha_e >0$. A real Landau singularity is a Landau singularity
 associated to a real Landau pattern. 
 
 \vskip .2cm
 
 (b)
 A {\em real Landau graph} is a  Landau graph $(\wt\Gamma, \wt\bphi)$
 (Defininiton \ref {def:landau-graph}) such that $\wt\bphi$ takes values in
 $V_\RR$ and  there exists a scattering diagram $\sigma = \bigl( (p_l), (q_e)\bigr)$
 of type $\wt\Gamma$ satisfying:
 \begin{itemize}
 \item[(1)] Each $p_l\in V$ lies in the real positive direction of the half-line
$\wt\bphi(l)$.

 \item[($2_\RR$)] For each $e=[v,w]\in \Ed(\Gamma)$ oriented from $v$ to $w$
 we have
 \[
 q_e\,\, \in\,\, \bigl( \RR_{>0}\cdot (\wt\bphi(w)-\wt\bphi(v))\bigr) \otimes
 \omega_{v,w}.
 \]
 Here $\omega_{v,w}\in \orr(e)$ is the basis vector corresponding to
 the orientation from $v$ to $w$. 
 \end{itemize}
 \end{defi}
 
 \noindent  For example, 
  Fig. \ref {fig:landau-graph} represents a real Landau graph,
  as the directions of the vectors in the local corolla satisfy ($2_\RR$). 
  
  \vskip .2cm
 
The restriction  to $\alpha_e>0$  in part (a) of Definition \ref {def:real-landau}
 and to
$\RR_{>0}$ in part (b)  is  justtified by
using the representation \eqref{eq:I(p)-Feynman}
which involves only non-negative $\alpha_i$. 
. It is this
 restriction that leads to appearance of convex geometry in the theory. 
 Proposition \ref {prop:Lsing=Lgrap} refines immediately to: 
 
  \begin{prop}
 A collection $p=(p_l)_{l\in L}$ represents a  real Landau singularity
 for $\wt\Gamma$
 if and only if it  comes from a real  Landau graph $(\wt\Gamma, \wt\phi)$, i.e., $p$ is 
 the external corolla of a  $V_\RR$-valued scattering diagram $\sigma$ 
  satisfying   the conditions of  Definition   \ref {def:real-landau} (b)
\qed
 
 \end{prop}
 
 In the sequel we will consider only real Landau patterns and singularities. 
 
 \paragraph{ ``Nearest'' singularities.} 
In this and the next paragraph we
  further make the following

 \begin{ass}
 At least one of the masses $m_e$ is not zero. 
\end{ass}

 With this assumption, we have: 
 
 \begin{prop}\label{prop:0-not-LS}
 The point  $p=0\in V^{L,0}_\RR$ is not a  real Landau singularity of $I(p)$. 

\end{prop}

\noindent{\sl Proof:} 
In the notation of \S \ref {sec:landau-complex}\ref {par:landau-geom}, the point
$0$ being a Landau singularity means that there is $\phi\in C^0(\Gamma, V)$
such that $\DD(\phi) = (\del \cdot \alpha^{-1} \cdot \delta) (\phi)=0$. 
For a real Landau singularity, $\phi\in C^0(\Gamma, V_\RR)$. 
 As $\delta = \del^*$,
this implies that
\[
0 = \bigl(\del \alpha^{-1} \delta \phi, \phi\bigr)\,=\,
(\alpha^{-1}\delta\phi, \delta\phi)\,=\, \sum_e \bigl( \alpha^{-1}_e (\delta\phi)_e, (\delta\phi)_e
\bigr)\,=\,\sum_e \alpha_e (q_e, q_e). 
\]
But $(q_e, q_e) = q_e^2$ must be equal to $m_e^2 \geq 0$, so the sum cannot be $0$
if at least one $m_e\neq 0$. . \qed

\vskip .2cm

We therefore have the maximal connected open set $U\subset V^{L,0}_\RR$
which contains $0$ and does not contain any  Landau singularities
satisfying our assumptions. The boundary $B=\del U$  of this set can be called the
{\em set of nearest (to $0$) Landau singularities} of $I(p)$. 
Consider the function
\[
\Phi(\alpha, p) \, = \, \Crv \bigl( F(\alpha, k)|_{\del^{-1}\eps_*(p)}\bigr).
\]
Although $F(\alpha, k)$ depends on $\alpha$ linearly, $\Phi(\alpha, p)$ does not. 
The following is a remark of  \cite[\S2]{landau1}. 

\begin{prop}
If $p^{(0)}\in B$ and $\alpha^{(0)}\in \RR_{>0}^{\Ed(\Gamma)}$ is the corresponding set of Landau
multipliers, then  the function $\alpha \mapsto \Phi(\alpha, p^{(0)})$  achieves a minimum at 
$\alpha^{(0)}$ with value $0$. 

\end{prop}

\noindent{\sl Proof:}  On the linear space $\del^{-1}(0)$,
the function $F(\alpha, k)$ achieves its critical value at $0$,
as it is the sum of a purely quadratic part and a constant.
This means    that  $\Phi(\alpha, 0) = \sum_e \alpha_e m_e^2$
which is strictly positive,  if all $\alpha_e >0$. This implies that
$\Phi(-, p): \alpha \mapsto \Phi(\alpha, p)$ is positive on $\RR_{>0}^{\Ed(\Gamma)}$ whenever $p\in U$. 
Indeed, otherwise the function $k\mapsto F(\alpha, k)|_{\del^{-1}\eps_*(p)}$ has, for 
some $\alpha\in
\RR_{>0}^{\Ed(\Gamma)}$, a critical point with critical value $0$
which means that $p$ is a Landau singularity which is not the case.
 Now, if $p^{(0)}$ actually
is a Landau singularity, then $\Phi(\alpha, p^{(0)})$ must be zero
for some  $\alpha\in \RR_{>0}^{\Ed(\Gamma)}$. 
If, moreover, $p^{(0)}\in B$, then for the nearby $p\in U$
we have $\Phi(-, p)>0$   on $\RR_{>0}^{\Ed(\Gamma)}$
and so for $\Phi(-, p^{(0)}$ to achieve $0$ on this set, it must achieve it
as a minimum. \qed

\paragraph{ Momenta in a subspace.} 
Let $W\subset V_\RR$ be any $\RR$-linear subspace.  It is often interesting to
study the case when the incoming momenta lie in $V$. 
 In particular, the {\em planar case},   when $\dim(W)=2$  was studied already
 in \cite{landau1}.   The following general statement is useful.

\begin{prop}
If $p=(p_l)_{l\in L}$ is a Landau singularity such that  $p_l\in W$ for all $l$,  then all the
intermediate momenta $q_e$ corresponding to  it,  lie in $W$ as well. 
\end{prop}

\noindent{\sl Proof:} Let $f: V_\RR\to \RR$ be any linear functional vanishing on $W$.
It is enough to show that $f(q_e)=0$ for any intermediate momentum. 
Let $\phi\in C^0(\Gamma, V_\RR)$ be the $0$-cochain giving rise to $p$. Then $\phi$
is not constant by Proposition \ref {prop:0-not-LS}. Let $v$ be such that $f(\phi(v))\geq 
f(\phi(w))$ for all $w\neq v$, and at least one of the inequalities is strict. 
Since $\Gamma$ is connected, we can assume
that $f(\phi(v)) > f(\phi(w))$ for some $w$ which is connected with $v$ by an edge. 
Since $f(\eps_*(p))=0$, we have
\[
0= f\bigl( \eps_*(p)_v\bigr) \,=\, f\bigl( (\del \alpha^{-1} \delta (\phi))_v\bigr) \,=\,
\sum_{e=[v,w]} \alpha_e^{-1} \bigl( f(\phi(w)) - f(\phi(v))\bigr).
\]
Here the summation is over all edges issuing from $e$.
If $f(\phi(v))\neq f(\phi(w))$, i.e., $f(q_e)\neq 0$ for at least one summand, the sum must be  strictly negative.  This proves the statement. \qed


\section{The Minkowski problem and momenta}\label{sec:minkowski}

Here and in the next section we relate real Landau singularities with convex geometry. 

\paragraph{Convex polytopes.}
  
  Let $E$ be a  real vector space of dimension $n$.
  By a {\em convex polytope} in $E$ we will mean a set $Q\subset E$
  which is the convex hull of a finite set of points. In particular, a convex
  polytope is bounded. For a convex polytope $Q$ we denote
  by $\Fc_k(Q)$ the set of $k$-dimensional faces of $Q$. Faces of
  codimension $1$ are called {\em facets}.

 \paragraph{Weighted normals as momenta.}  
 Put
  \[
  W\,=\,\Lambda^{n-1}(E) \, \= \, E^*\otimes\Lambda^n(E). 
  \]
   The space $E$ carries a tautological $(n-1)$-form $\Omega$ with
  values in $W$. This form corresponds to
   $\Id\in \Lambda^{n-1}(E^*) \otimes\Lambda^{n-1}(E)$. As $\Omega$ has
   constant coefficients, $d\Omega=0$. 
   
  \vskip .2cm
  
 Let  $Q\subset E$ be a convex polytope
  which we suppose to be of full dimension. 
  For each facet   $F\subset Q$ its
  {\em weighted normal}   is defined as
  \be
  \nu(F) = \nu_Q(F) \,=\, \int_F \Omega\,\, \in\,\,   W. 
  \ee
 Here we use the orientation on $F$ induced by the chosen orientation on $W$
 and the ``outward''  (with respect to $Q$) orientation of the normal to $F$.

   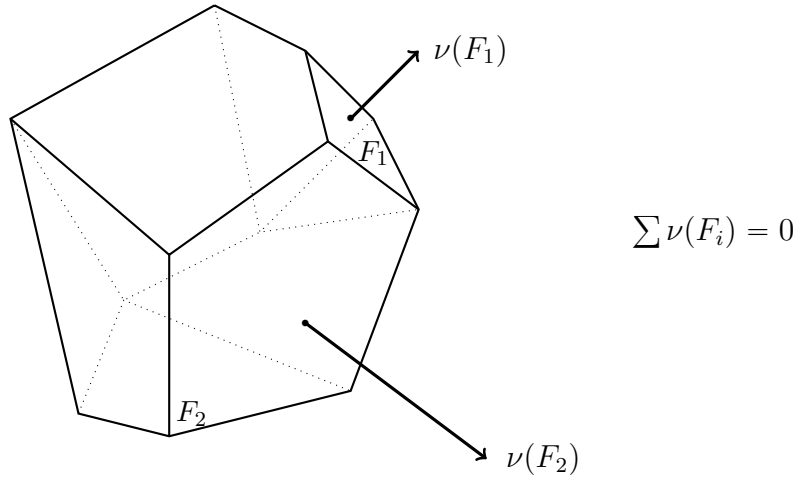
\begin{figure}[h]
 \centering
  \begin{tikzpicture}[scale=0.3]
  
  \draw [line width=0.8] (-2, 10) -- (2,8) --(5,5) -- (7,1) -- (4, -7) -- (-4, -9) -- (-8,-8) --  (-11,5) --(-2, 10); 
  \draw  [line width=0.8] (2,8) -- (3,4) -- (7,1); 
  \draw [line width=0.8]   (3,4) -- (-4, -1) -- (-11,5); 
  \draw  [line width=0.8]  (-4, -1) -- (-4, -9); 
  
  \node at (4,5) {\tiny$\bullet$}; 
  \draw[->, line width =1.2](4,5) -- (7,8); 
  
  \node at (2, -4){\tiny$\bullet$}; 
  \draw[->, line width =1.2] (2, -4) -- (10, -10); 
  
  \node at (5,3.5){\small $F_1$}; 
  \node at (9.3,8) {$\nu(F_1)$} ;

  \node at (-3, -8) {\small $F_2$}; 
  \node at (12.5, -10){$\nu(F_2)$}; 
  
  
  \node at (20, 0){$\sum \nu(F_i)=0$}; 
  
  \draw[dotted]   (-8, -8) -- (-6, -3) -- (4, -7) ; 
  \draw[dotted]   (-11,5) -- ((-6, -3) -- (0,0) -- (7,1); 
  \draw[dotted] (-2, 10) -- (0,0) -- (5,5); 
  
  \end{tikzpicture}
  \caption{Weighted normals in the classical sense.} 
  \label{fig:w-normals}
  \end{figure}
  
  \begin{prop} \label{prop:sum-normals=0}
  (a)  If $F$ is a common facet of two polytopes $Q, Q'$ situated
 on the opposite sides of $F$, then $\nu_Q(F)+\nu_{Q'}(F)=0$. 
 
 \vskip .2cm
 
 (b) For any polytope $Q$ we have 
  $
 \sum\limits_{F\in\Fc_{n-1}(Q)} \nu_Q(F) = 0. 
 $
 \end{prop}
 
 \noindent{\sl Proof:} (a) follows since the orientations of $F$ used in $\nu_Q(F)$ and 
 $\nu_{Q'}(F)$
 are opposite. Part (b) follows from the Stokes formula, since $d\Omega=0$. \qed
  
   \begin{ex}\label{ex:hodge-normals}
 Suppose that $W$ is equipped with a  Euclidean scalar product
 and an orientation. Then the Hodge $*$-operator
 identifies $W\= \Lambda^{n-1}(W)$.  After this identification,  $\nu_Q(F)$ 
 is realized  as the  vector in $W$  normal to $F$,
 facing outwards with length  equal to the area of $F$,
 see Fig. \ref {fig:w-normals}.  By Proposition \ref{prop:sum-normals=0}(b),
  the sum of such
 vectors is $0$. 
 
 \end{ex}

  \begin{ex} \label{ex:norm-polygon}
   Suppose that $n=2$, i.e., that $Q$ is a plane polygon with vertices $v_1, \cdots, v_N$
 numbered clockwise.  In this case $\Lambda^{n-1}(W)=W$ without the need
 for any additional structure. The $1$-form $\Omega$ is the Maurer-Cartan form
 of $W$ as a Lie group. It follows that
   for a side  
  $F = [v_k, v_{k+1}]$ of $Q$ the vector $\nu_Q(F)$ is equal to 
  the side vector  $v_{k+1}-v_k$.   Proposition
  \ref{prop:sum-normals=0}(b) in this case  expresses the ``closure of the polygon'' 
  $\sum (v_{k+1}-v_k)=0$. 

 If we equip $W$ with a Riemannian metric and an orientation,
  then the $*$-operator acts on $W$ as the rotation  by $90^\circ$   anticlockwise
  and we get the  weighted  (outer) normals to the sides as in 
  Example \ref {ex:hodge-normals}. 
  
   \end{ex}

  \paragraph{ The Minkowski problem.} 
  
  The converse problem of finding a polytope
  from the data of   its weighted normals (which must sum up to $0$) is known as  the
  {\em Minkowski problem} \cite{minkowski, gu-yau, klain}. The classic result of Minkowski
   \cite{minkowski, klain}
  says that this problem always has an essentially unique solution:
  
  \begin{thm}\label{thm:minkowski}
  Let $I$ be a finite set and
   $p=(p_i)_{i\in L}$ be a set  of distinct nonzero vectors in $W$
  such that:
  \begin{itemize}
  \item[(1)] The $p_i$ span $W$.
  
  \item[(2)] No two  of the $p_i$ are  positive multiples of each other. 
  
  \item [(3)] $\sum_{i\in I} p_i=0$.

  \end{itemize}
  Then there exists a unique, up to a  translation, convex polytope
  $Q=Q(p)\subset E$ and a bijection $a: I\to \Fc_{n-1}(Q)$ such that $p_i = \nu_Q (a(i))$ for all 
  $i\in I$.  \qed
  \end{thm} 
  
  More precisely, Minkowski's result is formulated for the concept of
 the $\nu_Q(F)$  understood in the classical sense (Example
  \ref {ex:hodge-normals}), in the presence of a Euclidean metric and
  orientation but since the definition  $\nu_Q(F)$ can be done invariantly,
  one can deduce the invariant formulation from the classical one. 
  
  \vskip .2cm
  
 Here, (1) and (2)  are just genericity conditions. If (1) is violated, then we
 can restrict the consideration to the subspace   in $W$ spanned by the $p_l$
 and redefine $E$ accordingly. If (2) is violated,  we can sum
 the $p_l$ in each positive proportionality class so the condition will hold for the
 new set. 
 
 \vskip .2cm
 
  Using  the  terminology of Definition \ref{def-scat-corolla}, we can say that 
 each  convex polytope $Q$ gives a (real) scattering corolla $\nu_Q$ 
 with values in $W$ labelled by $\Fc_{n-1}(Q)$ and each (real) scattering
  corolla satisfying the conditions (1) and (2) of Theorem 
  \ref {thm:minkowski} can be realized in this way.

  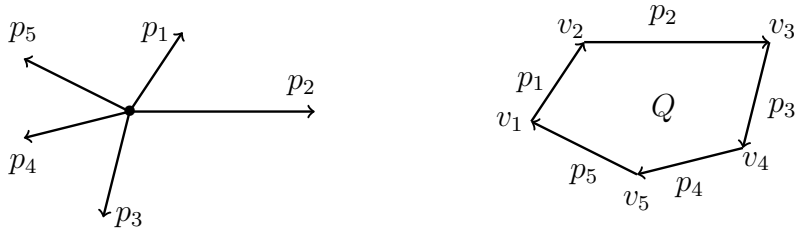
\begin{figure}[h]
 \centering
  \begin{tikzpicture}[scale=0.35]
  \node at (0,0){\small $\bullet$}; 
  \draw [->, line width=0.9]  (0,0)-- (2,3); 
   \draw [->, line width=0.9]  (0,0)-- (7,0); 
   \draw [->, line width=0.9]  (0,0)-- (-1, -4); 
    \draw [->, line width=0.9]  (0,0)-- (-4, -1); 
   \draw [->, line width=0.9]  (0,0)-- (-4,2);   
   
   \node at (1,3){$p_1$}; 
   \node at (6.5, 1) {$p_2$}; 
   \node at (0, -4){$p_3$}; 
   \node at (-4, -2){$p_4$}; 
  \node at  (-4, 3){$p_5$}; 
  \node at (10,0){};
  \end{tikzpicture}
  \begin{tikzpicture}[scale=0.35]
  \draw[->, line width=0.9] (0,0) -- (2,3); 
  \draw[->, line width=0.9] (2,3) -- (9,3);
  \draw[->, line width=0.9] (9,3) -- (8, -1);
  \draw[->, line width=0.9] (8, -1) -- (4, -2);
  \draw[->, line width=0.9] (4, -2) -- (0,0);
  
  \node at (-0.8,0){$v_1$}; 
  \node at (1.5, 3.5){$v_2$}; 
  \node at (9.5, 3.5){$v_3$}; 
  \node at (8.5, -1.5){$v_4$}; 
  \node at (4, -3){$v_5$}; 
  
  \node at (0, -4){}; 
  \node at (-4,0){}; 
  
  \node at (0, 1.5) {$p_1$}; 
  \node at (5,4) {$p_2$};
  \node at (9.5,0.5){$p_3$};  
  \node at (6, -2.5) {$p_4$}; 
  \node at  (2, -2){$p_5$}; 
  
  \node at (5,0.5){$Q$}; 
  \end{tikzpicture}
  \caption{ $2$D Minkowski problem : from a set of momenta to a convex polygon.}
  \label{fig:poly-momenta}
 \end{figure}

       \begin{ex}
   In the planar case considered in Example \ref  {ex:norm-polygon}, 
   the solution to the Minkowski problem is obvious. 
   That is, by rotating the $p_l$ by $90^\circ$ 
   clockwise, it suffices to construct a unique, up to translations, convex polygon 
 $Q=Q(p)\subset W$
 such that the $p_l$ are the vectors associated {\em to its sides}, read clockwise.
 That is, $Q$ has vertices $v_1,\cdots, v_N$, $N=|I|$,
  read clockwise
 and each  $v_{k+1}-v_k$ (with $k$ taken modulo $N$) is equal to a unique $p_{l_k}$ from the set, see Fig. \ref {fig:poly-momenta}. 
 To construct $Q$, we number  the $p_l$ clockwise,  
 according to their slopes
 as $p_1, \cdots, p_N$ (such numeration is unique up to
 a cyclic rotation) and then  put 
 \[
 v_1=0, \, v_2 = p_1, \, v_3 = p_1+p_2, \cdots, v_N = p_1 + p_2 +\cdots + p_{N-1} .
 \]
 The fact that $\sum p_l=0$ means that the polygon ``will close'',
 i.e., that the vector corresponding to the last side will be $p_N$. 
 
    \end{ex}


\section{Regular polyhedral subdivisions and real Landau singularities
 }\label{sec:RPS-and-RLS}

\paragraph{Polyhedral subdivisions and scattering diagrams.}  \label{par:polydec}

 Let us recall some formalism from  \cite{GKZ}.
 As before,  $E$  is a real vector space  of dimension $n$. Let us choose a translation
 invariant volume form on $E$,  thus trivializing $\Lambda^n(E)$ and
 identifying $W= \Lambda^{n-1}(E)$ with $E^*$. This also fixes an orientation
 on $E$. 

\vskip .2cm

Let $Q\subset E$ be a convex polytope which we assume of full dimension. 
   A
{\em polyhedral subdivision} of  $Q$ is a finite collection
$\Pc$  of convex polytopes (of full dimension)
such that $Q=\bigcup_{P\in \Pc} P $
and each intersection $P_1\cap P_2$, $P_1, P_2\in\Pc$
 is a proper face of them both (possibly empty). We denote by $\Fc_k(\Pc)$
 the set of $k$-dimensional faces of $\Pc$. Thus $\Fc_n(\Pc)$,
 the set of full-dimensional polytopes in $\Pc$,
  can be
 identified with $\Pc$ itself. We also denote $\Fc(\Pc) = \bigsqcup_k \Fc_k(\Pc)$
 to be the set of all the faces of $\Pc$. 
 
 \vskip .2cm
 
A {\em triangulation} is a polyhedral subdivision  into simplices. 
Given two polyhedral subdivisions $\Pc, \Pc'$, we say that
$\Pc'$ is a {\em refinement} of $\Pc$, if, for each polytope $P\in\Pc$,
the polytopes $P'$ from $\Pc'$ contained in $Q$,
form a polyhedral subdivision  of $P$. This will be expressed
by $\Pc'\preceq \Pc$. 

\vskip .2cm

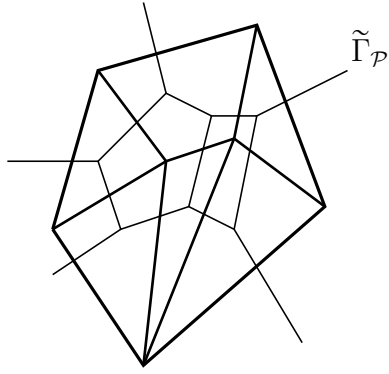
\begin{figure}[h]
 \centering
  \begin{tikzpicture}[scale=0.3]

  \draw[line width=1.2] (-5, -3) -- (-3,4) -- (4,6) -- (7, -2) -- (-1, -9) -- (-5, -3); 
  
  \draw[line width=0.6] (-3,0) -- (0,3) -- (2,2) --  (1, -2) -- (-2, -3) -- (-3,0); 
  \draw [line width=0.6] (2,2) -- (4,2) -- (3, -3) -- (1, -2); 
  
  \draw[line width=1] (0,0)-- (-5, -3); 
  \draw[line width=1] (0,0) -- (-3,4) ; 
  \draw[line width=1] (0,0) -- (3,1); 
  \draw[line width=1] (0,0) -- (-1, -9); 
  
  \draw[line width=1]  (3,1) -- (4,6); 
  \draw[line width=1] (3,1) -- (7, -2); 
  \draw[line width=1] (3,1) -- (-1, -9); 
  
  \draw [line width=0.6] (-7,0) -- (-3,0);
  \draw [line width=0.6] (-1, 7) -- (0,3); 
  \draw [line width=0.6] (8,4) -- (4,2);
  \draw [line width=0.6] (3,-3) -- (6,-8);
  \draw[line width=0.6]  (-2, -3) -- (-5, -5); 
  
  \node at (9,5) {$\wt\Gamma_\Pc$}; 
 
  \end{tikzpicture}
   
  \caption{A polyhedral subdivision  and its dual graph.} 
  \label{fig:polydec+graph}
  \end{figure}

 Given a polyhedral subdivision  $\Pc$ of $Q$,
  we have the {\em dual graph} $\wt\Gamma_\Pc$ of $\Pc$.
  By definition,
  \[
  \Vert(\wt\Gamma_\Pc) = \Fc_n(\Pc), \quad \Ed(\wt\Gamma_\Pc) = \Fc_{n-1}(\Pc)
  \]
  are the sets of the polytopes from $\Pc$ and of their  facets. 
  Thus  legs correspond to the facets  which lie in facets of $Q$ itself 
  (and so  enter into only one
 polytope $Q\in \Pc$) and other edges  (those of  the associated legless
 graph $\Gamma_\Pc$)
 correspond to internal facets
 (contained in exactly two of the polytopes), see Fig. \ref{fig:polydec+graph}.

  As explained in \S \ref {sec:feyn-setup}\ref {par:formal-momenta},  an element  
  of $C_1(\wt\Gamma_\Pc, W)$ is a datum $\sigma$   associating
 to any incident pair $(v,e)$ of a vertex and edge an element $\sigma(v,e)\in W$
 so that $\sigma(v,e) = -\sigma(v',e)$ whenever $e$ is an edge with two vertices
 $v,v'$. In our situation any such $(v,e)$ is given by a pair
   $P\supset F$ with $P\in \Fc_n(\Pc)$ and $F\in \Fc_{n-1}(P)$.
   So we can write $\sigma$ as $(\sigma(P,F))$. 
   
   \vskip .2cm

  Proposition  \ref{prop:sum-normals=0}  implies the following obvious but appealing
  
  \begin{prop}\label{prop:P-sc-diag}
  Let $\Pc$ be a polyhedral subdivision of $Q$. 
  Then the rule $\nu_\Pc$ associating to any pair $P\supset F$ as above
  the vector $\nu_P(F)$,  is a scattering diagram 
  (Definition \ref {def-scat-corolla}) with values in $W$ labelled by the graph $\wt\Gamma_\Pc$. 
  If $\Pc$ is such that the facets of $Q$ are not subdivided, then $\nu_\Pc$
  refines the scattering corolla $\nu_Q$. 
 \qed
   \end{prop}
   
   \begin{rem} Given any scattering diagram $(\wt \Gamma, \sigma\in C_1(\wt\Gamma, W))$
refining $\nu_Q$ (in particualar with $\Leg(\wt\Gamma)=\Fc_{n-1}(Q)$),
one can ask whether it comes from a subvidision $\Pc$ as above so that
$\wt\Gamma=\wt\Gamma_\Pc$. It is not true in general, one necessary (but not sufficient)
condition being this. Each vertex $v\in\wt\Gamma$ gives the local scattering corolla
$\sigma_v$ and so (assuming the genericity conditions of Theorem \ref{thm:minkowski})
a polytope $Q_{\sigma_v}$. Then we must have $\sum_v\Vol(Q_{\sigma_v}) = \Vol(Q)$. 

\end{rem}

 \paragraph{Convex functions and regular subdivisions.}

\vskip .2cm

Let $\Pc=(Q_i)$ be a polyhedral subdivision  of $Q$.
A continuous function $f: Q\to \RR$ will be called
{\em piecewise linear  (PL)} with respect to 
$\Pc$, if each restriction $f|_{Q_i}$ is an affine-linear function.
We denote by $\PL(\Pc)$ the real vector space formed by such
 functions. We also denote $\Aff(Q) = \Aff(W)$ the  space of
 affine-linear functions on $Q$ (equivalently, on the whole of $W$). 
 It acts on $\PL(\Pc)$ by addition. 
 
 \vskip .2cm
 
 A piecewise-linear function $f\in \PL(\Pc)$ will be called
 {\em convex}, if the region \hfill\break
  $\{(x,y) \in Q\times \RR \, | \,y\geq f(x)\}$
 is convex.  Such functions form a convex cone $\PL^+(\Pc)\subset
 \PL(\Pc)$ (not necessarily of full dimension). 
 
 \vskip .2cm

\begin{figure}[h]
 \centering
  \begin{tikzpicture}[scale=0.45]
  
  \draw (-5, -2) -- (0,4) -- (5, -2) -- (-5, -2); 
  \draw (-2, -1) -- (0,2) -- (2, -1) -- (-2, -1); 
  \draw(-5, -2) -- (-2, -1); 
  \draw (0,4) -- (0, 2); 
  \draw (5, -2) -- (2, -1); 
  \draw (-5, -2) -- (2, -1); 
  \draw (5, -2) -- (0,2);
  \draw (0,4) -- (-2, -1); 
  
  \end{tikzpicture}
  \caption{A non-regular triangulation. }
  \label{fig:non-reg-tr}
  \end{figure}
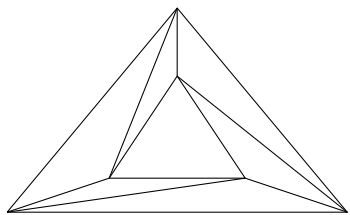

 We say that $f\in \PL^+(\Pc)$ is {\em strictly convex}, if $f$ actually
 breaks on each intermediate facet of $\Pc$, i.e., on each
 facet which is common to two polytopes in $\Pc$. 
 This means that it lies in the interior of $\PL^+(\Pc)$. 
 A subdivision  admitting a strictly convex PL-function is called
 {\em regular}.  Put differently, $\Pc$ is regular, if
 $\PL^+(\Pc)\subset \PL(\Pc)$ is of full dimension.  
  A classical example of a non-regular
 subdivision (a plane triangulation) is in Fig. \ref {fig:non-reg-tr}. 
 
  \paragraph{Regular subdivisions and real Landau singularities.} 
  
  Let us now assume that $E$ is equipped with a real scalar product
   $(-,-)$
  as in \S  \ref {sec:real-gen}
 \ref {par:alpha>0}. We further assume that the chosen volume form
 is compatible with the volume defined by the scalar product,
 so the identification $W=\Lambda^{n-1}(E) \= E^*$ preserves the
 product.

  \begin{thm}\label{thm:polyhedral}
  Let $\Pc = (Q_i)_{i\in I}$ be a polyhedral subdivision of $Q$ with dual graph $\wt\Gamma$.
  For each edge $e\in\wt\Gamma$ let $F_e$ be the corresponding facet of one or two
  polytopes, $P$ or $(P,P')$  from $\Pc$. We assume that  the scalar 
  product $(\nu_P(F_e), \nu_p(F_e)) \geq 0$  for all $e$ and let $m_e$ be its non-negative square root. 
  If $\Pc$ is regular, then the scattering diagram associated to $\Pc$ by
  Proposition \ref {prop:P-sc-diag}, is a real Landau singularity with the mass of a
  (virtual) particle associated to $e$ being  $m_e$. 
  \end{thm}
  
  Thus in  the Eucliden case, 
  $m_e$ is just the $(n-1)$-dimensional volume of $F_e$. 
  
  \vskip .2cm
  
  \noindent{\sl Proof:} We follow the procedure of \S \ref {sec:landau-complex} \ref {par:landau-geom} 
  on geometric construction of Landau singularities, restricting our attention to real data.

   Let $f: Q\to\RR$ be a strictly convex $\Pc$-piecewise-linear function. 
  For each polytope $P\in \Pc$  let $\phi(P) = d(f|_P)$ be the linear part of $f$ on $P$. It is
  an element of $E^*$ which we identified  with $W$. This gives a spatial realization
  $\phi$ of $\Gamma$ in $W$.  An oriented  edge $e=[v,v']$ of $\Gamma$ corresponds to an ordered pair
  $(P, P')$ of two polytopes from $\Pc$ with a common facet $F$. 
  The difference $(\delta\phi)_e$ is the {\em break} of $f$ on the facet $F$. 
  The fact that $f$ is continuous means that this difference actually vanishes on $F$. 
   Linear forms vanishing on (the linear subspace associated to) $F$ form a $1$-dimensional
   subspace $F^\perp\subset E^*=W$. This subspace is spanned by the
   weighted normal $\nu_P(F)$. The fact that $f$ is strictly convex now means that
   $(\delta\phi)_e$ is a positive multiple of $\nu_P(F)$, i.e., there is $\alpha_e>0$ such that 
   $(\delta\phi)e = \alpha_e\cdot \nu_P(F)$.  In other words, we obtained
   a collection of positive Landau multipliers thus producing a real Landau singularity. 
   \qed
 
 \vskip .2cm
 We will refer to Landau singularities given by Theorem
 \ref {thm:polyhedral} as {\em polyhedral}. 
 
 \paragraph{The Landau graph of a polyhedral singularity.} 
 \label{par:LG-polyh}
  For future reference
 let us describe the real Landau graph (Definition \ref {def:real-landau})
 corresponding to a polyhedral Landau singularity as an explicit
 embedding $\wt\bphi: \wt\Gamma_\Pc \to W^*$ into the dual space. 
 
 \vskip .2cm
 
 Let a convex function $f: Q\to \RR$ from $\PL^+(\Pc)$ be given and 
 let $l\in W^*$. Put
 \[
 h_l \,=\, \max \,\bigl\{ h\in \RR \,| \, l(x) + h \, \leq \, f(x), \,  \,\,\forall \, x\in Q\bigr\},    
 \]

    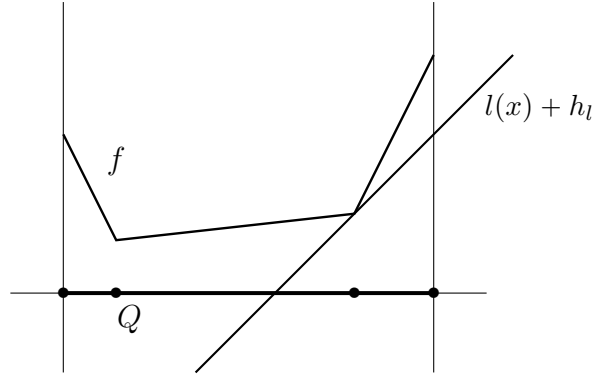
\begin{figure}[h]
 \centering
  \begin{tikzpicture}[scale=0.35]
  
  \draw[line width=1.5] (0,0) -- (14,0); 
  \draw (-2,0) -- (16, 0); 
  \draw (0, -3) -- (0,11); 
  \draw (14, -3) -- (14, 11); 
  \draw [line width = 0.9] (0,6) -- (2,2) -- (11,3) -- (14,9); 
  
  \draw  [line width = 0.8] (5, -3) -- (17,9); 
  
  \node at (0,0){\small$\bullet$}; 
    \node at (2,0){\small$\bullet$}; 
      \node at (11,0){\small$\bullet$}; 
        \node at (14,0){\small$\bullet$}; 
        
   \node at (18,7){\small $l(x) +h_l$};  
   
   \node at (2.5, -1){$Q$};  
   \node at (2,5){$f$};    
  
  \end{tikzpicture}
  \caption{The constant $h_l$. }
  \label{fig:h-l}
  \end{figure}

  see Fig. \ref {fig:h-l}. 
 Then the set
 \[
 F_l = \bigl\{ x\in Q\,\bigl| \, l(x)+h_l = f(x)\bigr\}
 \]
is a face of (some polytope from) $\Pc$. So for each $F\in\Fc(\Pc)$
we put
\[
\wc F \,=\, \bigl\{ l\in W^*\, | \,\, F_l \supseteq F\bigr\}. 
\]
  This is a  closed convex polyhedron (possibly unbounded)  in $W^*$ 
  and the $\wc F, F\in\Fc(\Pc)$,  form a  decomposition  of $W^*$ into possibly unbounded
  polyhedra  which we denote $\Sc_f$. 
  
  \begin{prop}\label{prop:LG-polyhedral}
  The $1$-skeleton $\wt\Gamma_f := \Sk_{\leq 1}(\Sc_f)\subset W^*$ is a real Landau graph 
  identified  with $\wt\Gamma_\Pc$. This graph produces the polyhedral Landau
  singularity corresponding to $(\Pc, f)$.    In particular,  Landau graphs
  producing polyhedral Landau singularities, do not have self-intersections. 
  \end{prop}
  
  \noindent{\sl Proof:} Indeed, vertices of $\wt\Gamma_f$ correspond to the maximal
  faces of $\Pc$, i.e., to the polytopes $P\in \Pc$ of dimension $n$. More precisely,
  these vertices are the differentials $d(f|_P)\in W^*$, $P\in \Fc_n(\Pc)$. 
  If two such polytopes $P, P'$ have a common facet $F$, then $\wt F$
  is the straight interval joining $d(f|_P)$ and $d(f|_{P'})$. If $F$ is an external facet
  of $P$, i.e., it is also a facet of $Q$, then $\wt F$
  is the infinite half-line 
  $d(f|_P) + \RR_+ \cdot d(g)$, where $g$ is any affine-linear function on $W$
  such that $g|_F=0$ and $g|_P\leq 0$.  This identifies $\wt\Gamma_f$ with
  $\wt\Gamma_\Pc$, i.e., defines a piecewise-affine
   embedding $\wt\bphi: \wt\Gamma_\Pc\to W^*$
  with image $\wt\Gamma_f$. The fact that $(\wt\Gamma_\Pc, \wt\bphi)$
  is a real Landau graph representing our singularity, follows directly
  from Definition \ref {def:real-landau} and from the proof of Theorem
  \ref{thm:polyhedral}. \qed

\paragraph{Landau multipliers and the Aleksandrov problem.}
Landau multipliers appearing in Theorem \ref {thm:polyhedral}
have an interpretation in terms of another classical problem
of convex geometry which can be called the {\em Aleksandrov problem}. 
For simplicity,  assume that the scalar product  on $E$ is Euclidean. 

\vskip .2cm

Let $Q\subset E$ be a convex polytope as before. 
A strictly
convex PL-function $f: Q\to \RR$  with respect to some polyhedral subdivision can be written as
\be\label{eq:f-as-min}
f(x)  \,=\, \max_{i=1}^N  \,\,  (l_i(x) + h_i),
\ee
 where $l_i: E\to\RR$, $i=1,\cdots, N$,   are linear functionals 
 and $h_i\in \RR$
 are constants.  Conversely, given  any $l_i, h_i$, we have a
 regular polyhedral subdivision $\Pc$ of $Q$ consisting
 of the polytopes
  \[
 P_i \, =\, \bigl\{ x\in Q\,\bigl| \, d_x f = l_i\bigr\} 
 \]
 (we count only those of $P_i$ that are of full dimension)
 for which $f$ defined as above is strictly convex. 
 Now, suppose the $l_i$ are fixed and we are
 given posittive numbers
 $\bv_i > 0$ such that $\sum_{i=1}^N \bv_i = \Vol(Q)$.   The Alexandrov problem consists
 in finding the  $h_i$ such that the  polytope $P_i$ has volume $\bv_i$
 for each $i$. It can be seen as a version of the Minkowski problem
 but for unbounded polyhedra of the form 
 $\{(x,y) \in Q\times \RR \, | \,y\geq f(x)\}$. 
 The fundamental theorem of Aleksandrov
 \cite{aleksandrov, gu-yau} can be formulated as:
 
 \begin{thm}
 The Aleksandrov problem always has a unique solution. \qed
 \end{thm}
 
 Differential properties of  the correspondence $(h_i) \mapsto (\bv_i := \Vol (P_i))$ is
  were studied
  by Gu, Luo, Sun and Yau
 \cite{gu-yau}.  
 To formulate them, let us assume that $E$ and therefore $E^*$
 is equipped with a Euclidean scalar product. 
  For $l\in E^*$ we denote by $\|l\|$ its Euclidean
length. Then \cite[Prop.2.4] {gu-yau} gives:
 
 \begin{prop}
 For $i\neq j$ we have
 \[
 {\del \bv_i \over \del h_j} \,=\, {A_{ij} \over \|l_i-l_j\|}.  
 \]
Here $A_{ij}$ is the area ($(n-1)$-dimensional volume)
 of the common face $F_{ij}$ between $P_i$ and $P_j$,
 in particular, $A_{ij}=0$ if $P_i$ and $P_j$ do not
 have a common facet. \qed
  \end{prop}
  
  In the situation of  Theorem \ref {thm:polyhedral}, we can reformulate
  the above proposition as:

 \begin{cor}
 Let $\Pc$ be a regular polyhedral subdivision of $Q$
 with dual graph $\wt\Gamma_\Pc$. Let us number the polytopes
 from $\Pc$ as $P_1,\cdots, P_N$. Consider the polyhedral
 Landau singularity corresponding to a strictly convex 
  $f\in \PL^+(\Pc)$ which we represent  in the form \eqref{eq:f-as-min}
Then for any  edge $e$ of $\Gamma_\Pc$ corresponding to
 a pair of  adjacent polytopes $P_i, P_j$, 
 the Landau multiplier $\alpha_e$ is equal to $\del \bv_i/\del h_j$. 
 
 \end{cor}
 
 \noindent{\sl Proof:}  We inspect the proof of Theorem \ref {thm:polyhedral} in our  modified notation. 
  Denoting by $v_i$ the vertex of $\Gamma_\Pc$
 corresponding to $P_i$, the $0$-cochain $\phi$ giving the spatial
 realization of $\Gamma_\Pc$ has the form $\phi(v_i) = l_i$.
 So for  an edge $e =[v_i, v_j]$  we have $(\delta\phi)_e = l_i-l_j$.
If $F$ is the common face of $P_i$ and $P_j$, then $A_{ij}$
is the length of $\nu_{P_i}(F)$. The Landau multiplier $\alpha_e$
is the coefficient of proportionality between the collinear covectors
$\nu_{P_i}(F)$ and  $(\delta\phi)_e$, i.e., the ratio
$A_{ij}/(\|l_i-l_j\|) = \del \bv_i/\del h_j$. \qed


\section{The planar case  }\label{sec:planar}

Here we consider the case when the momenta lie in a $2$-dimensional real subspace $W\subset V_\RR$.
This situation was studied geometrically in \cite{landau1, okun} under the name ``method of schemes''. 
We present a more detailed formal  statement behind this method. 
In this case it is possible to characterize the class of polyhedral
Landau singularities more explicitly. As $\dim(W)=2$, we use the term    {\em polygonal} instead of ``polyhedral''
throughout this section.  We also fix an orientation of $W$.

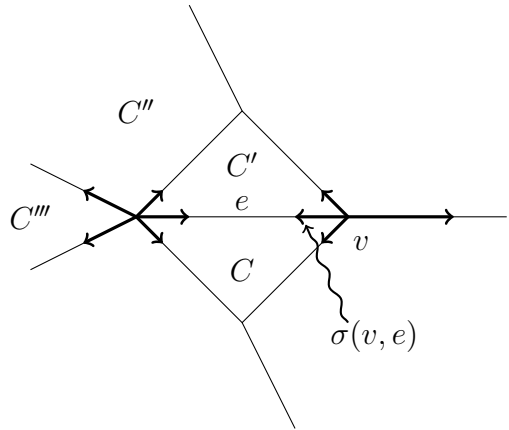
\begin{figure}[h]
 \centering
  \begin{tikzpicture}[scale=0.35]
 
 \draw (0,0) -- (4,4) -- (8,0) -- (4, -4) -- (0,0); 
 \draw (0,0) -- (8,0); 
 \draw (-4,2) -- (0,0) -- (-4, -2); 
  \draw (8,0) -- (14,0); 
   \draw (4, -4) -- (6, -8);
    \draw(4,4) -- (2,8); 
    
 \node at (4, -2){$C$};   
  \node at (4, 2){$C'$};  
  \node at (0, 4){$C''$};  
   \node at ( -4,0) {$C'''$}; 
   
   \draw [->, line width = 1.2] (0,0) -- (1,1);
    \draw [->, line width = 1.2] (0,0) -- (1, -1);
    \draw [->, line width = 1.2]  (0,0) -- (2,0); 
       \draw [->, line width = 1.2] (0,0) -- (-2,1); 
     \draw [->, line width = 1.2] (0,0) -- (-2, -1);

      \draw [->, line width = 1.2] (8,0) -- (6,0);   
    \draw [->, line width = 1.2] (8,0) -- (7,1); 
    \draw [->, line width = 1.2] (8,0) -- (7, -1); 
     \draw [->, line width = 1.2] (8,0) -- (12,0); 
     
     \node at (8.5, -1){$v$}; 
    \node at (4, 0.5){$e$};  
     
     \draw[decorate,
      decoration={snake, amplitude=0.4mm, segment length=4mm},->, 
     line width = 0.8] (8, -4) to 
     (6.4, -0.3);

     \node at (9, -4.5){$\sigma(v,e)$}; 
       
 \end{tikzpicture}
 \caption{A planar Landau graph giving a  cell decomposition of the plane,
 with two corollas of its scattering diagram depicted. }
 \label{fig:planar-LG}
 \end{figure}

\paragraph{Planar Landau graphs and polygonal subdivisions.} 
Let us call a real Landau graph $(\wt\Gamma, \wt\bphi: \wt\Gamma\to W)$
 (Definition \ref {def:real-landau}) {\em planar}, if $\wt\bphi$ is an embedding,
 i.e., the image  $\wt\bphi(\wt\Gamma)$ does not have self-intersections.
We can then identify $\wt\Gamma$ with its image and think of it
 as a subset of $W$.  A non-planar Landau graph with values in
 a $2$-dimensional space is depicted
 in  Fig. \ref{fig:landau-graph}, a planar one in Fig. \ref {fig:planar-LG}. 
 
 \begin{prop}\label{prop:polyg=planar}
 A real Landau singularity with momenta in $W$ is polygonal if and only if
 the corresponding Landau graph is planar. 
 \end{prop}
 
 \noindent{\sl Proof:} The ``only if'' part is a particular case of
 Proposition \ref{prop:LG-polyhedral}. To prove the ``if'' part, suppose that
 $\wt\Gamma\subset W$ is a planar Landau graph.  
 
 \vskip .2cm As $\wt\Gamma$ is connected and has at least $2$ legs,
  $W\-\wt\Gamma$ is the union of contractible regions, i.e., of $2$-dimensional
 cells. Together with $\wt\Gamma$, this forms a cell decomposition, denote
 it $\Sc$, of $W$  (with some cells unbounded), see Fig. \ref {fig:planar-LG}. 
 
\noindent  Denote by $\Sc_k$
 the set of $k$-cells of $\Sc$, so $\Sc_0=\Vert(\wt\Gamma)$ and
 $\Sc_1 = \Ed(\wt\Gamma)$. Let $C_\bullet(\Sc,W)$ be the cellular chain
 complex of $\Sc$ with coefficients in $W$:
 \[
 C_2(\Sc, W) = \bigoplus_{C\in\Sc_2} W \buildrel \wt\del\over\lra C_1(\wt\Gamma, W)
 \buildrel \wt\del\over\lra C_0(\wt\Gamma, W). 
 \]
It calculates $H_\bullet^\BM(W,W)$, the Borel-Moore homology of
$W$ as a topological space with coefficients in $W$ as an abelian group,
which has the form
\[
H_k^\BM(W,W) = H_k^\BM(W,\RR)\otimes_\RR W
=\begin{cases}
W, & \text{ if } k=2;
\\ 0, & \text {otherwise}. 
\end{cases}
\]
Therefore the scattering diagram $\sigma$ of our Landau singularity, being a
cycle in $C_1(\wt\Gamma, W)$, must be a boundary in $C_\bullet(\Sc, W)$.
That is, $\sigma = \wt\del(\tau)$ where $\tau = (\tau(C))_{C\in \Sc_2}$ is a collection
of vectors of $W$ associated to the $2$-cells. In other words, whenever $C$ and $C'$
are adjacent $2$-cells with an edge $e$ between them so that crossing
from $C$ to $C'$ leaves a vertex $v\in e$ on the right (with respect to our
orientation of $W$), we have
\[
\tau(C) - \tau (C') = \sigma(v,e)
\quad \text{(notation from  \S \ref {sec:feyn-setup} \ref {par:formal-momenta},
see also Fig. \ref {fig:planar-LG})}. 
\]
By going around the unbounded cells $C$, this implies that the $\tau(C)$ for
such $C$ form the vertices of a convex polygon $Q=Q(p)\subset W^*$ associated, via the
Minkowski problem, to the scattering corolla $p=(p_l)_{l\in L}$
formed by the incoming moments of $\sigma$.

\vskip .2cm

Similarly, by going around the cells $C$
containing some fixed vertex $v\in\wt\Gamma$ in their closure,
this implies that the $\tau(C)$ for such $C$ form the vertices of
a convex polygon $P_v := Q(\sigma_v)\subset W^*$ associated, via the Minkowski problem,
to the local  scattering orolla $\sigma_v$ of $\sigma$. 

\vskip .2cm 

Further, this means that the $P_v$ form a polygonal
subdiviion $\Pc$ of $Q$ with the dual graph $\wt\Gamma_\Pc$ identified with
$\wt\Gamma$ as an abstract graph. 

\vskip .2cm 

Finally, it remains to construct a strictly convex PL-function $f: Q\to \RR$
so that $\Sc=\Sc_f$ is the cell decomposition associated to $f$ as in 
\S \ref {sec:RPS-and-RLS} \ref {par:LG-polyh}. We know the slopes
that $f$ must have: they are the linear functions $l_v: W^*\to\RR$
associated to the vertices $v$ of $\Gamma$ which are points of $W$.
So there should be some $h_v\in \RR$ such that 
\be\label{eq:f_h-2dim}
f(x) = l_v(x) + h_v, \quad\text{when} \quad x\in P_v.
\ee
The obvious necessary condition on the $h_v$ is that $f$ should
be continuous, i.e., that whenever $P_v$ and $P_w$ have a common side $F$, 
the definitions \eqref {eq:f_h-2dim} should give the same function on $F$.
Now, such $F=F_e$ corresponds to an edge $e=[v,w]$ of $\wt\Gamma$. 
By construction of the faces (sides) of the $P_v$ via the Minkowski problem,
the linear parts $l_v, l_w$  coincide on the $1$-dimensional
linear subspace (line through $0$) parallel to $F_e$. This means that
for $x$ lying on  $F_e$ itself we have $l_v(x) -l_w(x) = c_{v,w}$
is a constant. 

\vskip .2cm

This gives the condition of continuity along $F_e$
as $h_v-h_w=c_{v,w}$. Note that $c_{v,w}=-c_{w,v}$, so
the system $c= (c_{v,w})$, a $1$-cochain (cocycle)
of $\Gamma$ (the legless graph obtained from $\wt\Gamma$)
with coefficients
in $\RR$.  To show existence of $(h_v)$ we must show
that $c$ is  a coboundary.  This means that the pairing of $c$ with
any $1$-cycle in $\Gamma$ vanishes. As $\Gamma$ is embedded into the
plane, a basis of cycles is provided by edge paths along the boundaries
of  bounded $2$-cells $C\in\Sc_2$. So we need to show that the sum of the
$c_{v,w}$ along the boundary of any such $C$ is $0$. That is, if we order
the vertices along the boundary as $v_1, \cdots, v_m$ clockwise
(with $v_{m+1}=v_1$), then
\[
\sum_{i=1}^m c_{v_i, v_{i+1}} = 0. 
\]
Now, $C$ corresponds to a vertex $\tau(C)$ of one of the polygons from $\Pc$
and the edges $e_i = [v_i, v_{i+1}]$ correspond to the sides $F_i = F_{e_i}$
of the polygons from $\Pc$, all passing through $\tau(C)$. As 
\[
c_{v_i, v_{i+1}}= l_{v_i}(x) - l_{v_{i+1}}(x), \quad \forall \,\, x\in F_i 
\]
we can specialize to $x=\tau(C)$ common to all the $F_i$, and obtain 
\[
\sum c_{v_i, v_{i+1}} \,=\, \sum \bigl(  l_{v_i}(\tau(C)) - l_{v_{i+1}}(\tau(C))\bigr) \,=\, 0.
\]
This shows that  $c$ is indeed a coboundary and so 
there is a unique, up to an additive constant, system
$h=(h_v)$ such that \eqref {eq:f_h-2dim} defines a 
{\em continuous} PL-function $f:Q\to\RR$.  It remains to show that $f$ is
convex, i.e., its breaks along the intermediate faces of $\Pc$ are positive
(``going up'').
But this follows from the definition of the polygons $P_v$ via the
Minkowski problem. \qed

\paragraph{Relation to plane webs of Gaiotto-Moore-Witten.} The concept of
planar Landau graphs is  essentially  equivalent to the more recent (2015)
notion of {\em plane webs} introduced by Gaiotto, Moore and Witten
\cite{GMW-big, GMW-small}. The motivation of  \cite{GMW-big, GMW-small}
was  quite different: 
study of scattering
of solitons and instantons in supersymmetric $2$-dimensional 
 Landau-Ginzburg models and more general $\Nc=(2,2)$ supersymmetric theories. 
 
 \vskip .2cm
 
 According to \cite{GMW-big, GMW-small}, see also \cite[Def.13.1]{KKS},  a plane web labelled by
 a finite set $\VV\subset \CC = \RR^2$  a datum of: 
 \begin{itemize}
 \item [(1)] 
  A planar graph $\wt\Gamma\subset \CC$ with edges represented by straight
  intervals or 
  half-lines.
  
  \item [(2)] A marking  assigning to any $2$-cell $C\subset \CC\-\wt\Gamma$
  an element $z_C\in \VV$, so that:
  
  \item[(3)] If an edge $e\in \Ed(\wt\Gamma)$ is oriented so that a $2$-cell $C$ is 
  on its left and $C'$ is on its right, then $e$ is parallel (and having the same  to the vector
  $z_C-z_{C'}$ and has the same direction. 
 
 \end{itemize}
 
 \noindent By viewing $\CC$ as a $2$-dimensional real vector space $W$ and
 viewing $z_C$ as element of $W^*$, we see that the concept is identical to
 the planar version of Definition \ref {def:real-landau} of real Landau graphs.
 The  scattering corolla  $\sigma_v$ corresponding to a vertex $v$ of
 $\wt\Gamma$ is formed by the vectors $z_C-z_{C'}$ for adjacent pairs of
 $2$-cells $(C, C')$ going clockwise around $v$. The Landau mutliplier
 $\alpha_e$ is the proportionality coefficient between the vector $\vec{e} = w-v$
 and $z_C-z_{C'}$ for an  oriented edge $e$ as in the condition (3) above. 
 
 \vskip .2cm
 
 The essential identity of the language of plane webs and that of
 regular polygonal subdivision has been pointed out earlier in \cite{KKS}. 
 The proof of Proposition \ref {prop:polyg=planar} given above is just a more
 detailed version of the argument in \cite[\S13]{KKS} adjusted to the terminology
 of Landau graphs.  
 
  \vskip .2cm
  
  It is remarkable that the same  mathematical structure appears in
  such seemingly different physical problems.


	\vskip 1cm

Kavli IPMU (WPI), UTIAS, The University of Tokyo,
5-1-5 Kashiwanoha, Kashiwa, Chiba, 277-8583 Japan. 

Email: {\tt mikhail.kapranov@protonmail.com}

\end{document}